\documentclass[a4paper]{amsart}
\usepackage{amscd, amssymb, color, booktabs, chemarrow}
\usepackage{myown}

\usepackage{harpoon}

\title{Degree Differences in the Eta Correspondences}
\author{Shu-Yen Pan}
\address{Department of Mathematics,
National Tsing Hua University, Hsinchu 300, Taiwan}
\email{sypan@math.nthu.edu.tw}


\date{\today}

\begin{document}

\begin{abstract}
A sub-relation of the $\Theta$-correspondence called the \emph{$\eta$-correspondence} is defined by
Gurevich-Howe for a finite reductive dual pair in stable range.
In this paper we propose an extension of the correspondence to general finite reductive dual pairs.
Then we determine the domain the correspondence and prove a formula on the difference of degrees in $q$ of
two irreducible characters paired by the correspondence.
\end{abstract}

\maketitle
\tableofcontents

\section{Introduction}

\subsection{}
Let $(\bfG,\bfG')$ be a reductive dual pair over a finite field $\bfF_q$ of $q$ elements
where $q$ is a power of an odd prime.
By restricting the character of the Weil representation to the dual pair,
we obtain a decomposition
\[
\omega_{\bfG,\bfG'}=\sum_{\rho\in\cale(G),\ \rho'\in\cale(G')}m_{\rho,\rho'}\rho\otimes\rho'
\]
where $\cale(G)$ denotes the set of irreducible character of the finite group $G$
of rational points of $\bfG$.
The subset
\[
\Theta_{\bfG,\bfG'}=\{\,(\rho,\rho')\mid m_{\rho,\rho'}\neq 0\,\}\subset\cale(G)\times\cale(G')
\]
gives a relation between $\cale(G)$ and $\cale(G')$, and is called the \emph{$\Theta$-correspondence}
for the dual pair $(\bfG,\bfG')$.
We say that $\rho'\in\cale(G')$ \emph{occurs} in the $\Theta$-correspondence for $(\bfG,\bfG')$
if the set
\[
\Theta(\rho')=\{\,\rho\in\cale(G)\mid m_{\rho,\rho'}\neq 0\,\}
\]
is not empty.

For a non-negative integer $k$, let $\bfG_k$ denote one of the following types of classical groups:
a general linear group $\GL_k$, a unitary group $\rmU_k$, a symplectic group $\Sp_k$ (only when $k$ is even),
or an orthogonal group $\rmO^\epsilon_k$ where $\epsilon=+$ or $-$.
For $\rho'\in\cale(G')$, the \emph{$\Theta$-rank} of $\rho'$,
denoted by $\Theta\text{\rm -rk}(\rho')$,
is defined to be the smallest number $k$ such that $\rho'\chi'$ occurs in the $\Theta$-correspondence
for some dual pair $(\bfG_k,\bfG')$ and for some linear character $\chi'$ of $G'$
(\cf.\ (\ref{0310}), (\ref{0311}) and (\ref{0405})).

For a dual pair $(\bfG_k,\bfG'_n)$ is stable range,
Gurevich and Howe prove in \cite{gurevich-howe} and \cite{gurevich-howe-rank} that for $\rho\in\cale(G_k)$
there is a unique $\rho'\in\Theta(\rho)$ such that $\Theta\text{\rm -rk}(\rho')=k$.
The character $\rho'$ is denoted by $\eta(\rho)$,
and the mapping $\rho\mapsto\eta(\rho)$ from $\cale(G_k)\rightarrow\cale(G'_n)$ is injective
and is called the \emph{$\eta$-correspondence} for the dual pair $(\bfG_k,\bfG'_n)$ in stable range.

In \cite{pan-eta} and \cite{pan-eta-unitary},
a one-to-one sub-correspondence $\underline\theta$ of the $\Theta$-correspondence between
$\cale(G_k)$ and $\cale(G'_n)$ for a general dual pair $(\bfG_k,\bfG'_n)$ is defined (\cf.~(\ref{0301})),
and it is shown that $\eta=\underline\theta$ if dual pair $(\bfG_k,\bfG'_n)$ is in stable range.
Therefore, for general dual pair $(\bfG_k,\bfG'_n)$ (i.e., not necessarily in stable range),
it might be reasonable to define $\eta(\rho)=\underline\theta(\rho)$ on those $\rho\in\cale(G_k)$
such that $\underline\theta(\rho)$ is defined and $\Theta\text{\rm -rk}(\rho)=k$.

For a non-negative integer $\ell$,
a condition on irreducible characters $\rho\in\cale(G)$ called
\emph{$\ell$-admissible for a dual pair $(\bfG,\bfG')$} is defined in
Subsections~\ref{0305}, \ref{0415} and \ref{0510}.
Then we have our first main result which gives complete description of the domain of the
``$\eta$-correspondence'' for a general dual pair:

\begin{thm}\label{0101}
Suppose that $0\leq k\leq n$.
\begin{enumerate}
\item[(i)] If $\rho\in\cale(G_k)$ is $(n-k)$-admissible for the dual pair $(\bfG_k,\bfG'_n)$,
then $\underline\theta(\rho)$ is defined and $\Theta\text{\rm -rk}(\underline\theta(\rho))=k$.

\item[(ii)] If $\rho'\in\cale(G'_n)$ is of\/ $\Theta$-rank $k$,
then $\rho'=\underline\theta(\rho)\chi'$ for
\begin{itemize}
\item some linear character $\chi'$ of $G'_n$; and

\item some group $\bfG_k$ and an irreducible character $\rho\in\cale(G_k)$ which is $(n-k)$-admissible for
the dual pair $(\bfG_k,\bfG'_n)$.
\end{itemize}
\end{enumerate}
\end{thm}

\begin{rem}
\begin{enumerate}
\item A condition on an irreducible character $\rho'\in\cale(\GL_n(q))$ of $\Theta$-rank $k$
is given in \cite{gurevich-howe-rank} theorem 9.2.3.

\item A similar condition on an irreducible character $\rho'\in\cale(\rmU_n(q))$ of $\Theta$-rank $k$
can be found in \cite{GLT01} theorem 3.9.
\end{enumerate}
\end{rem}

\subsection{}
For $\rho\in\cale(G_k)$,
it is known that $\rho(1)$ is a polynomial in $q$.
Let $\deg_q(\rho)$ denote the degree of this polynomial.
It turns out that we have a very elegant formula on the difference between
$\deg_q(\rho)$ and $\deg_q(\underline\theta(\rho))$ when $\Theta\text{\rm -rk}(\underline\theta(\rho))=k$:

\begin{thm}\label{0501}
Consider a dual pair $(\bfG_k,\bfG'_n)$ where $k\leq n$.
If $\rho\in\cale(G_k)$ is $(n-k)$-admissible for $(\bfG_k,\bfG'_n)$,
then
\[
\deg_q(\underline\theta(\rho))=\deg_q(\rho)+
\begin{cases}
\frac{1}{2}k(n-k+1), & \text{if\/ $(\bfG_k,\bfG'_n)=(\rmO^\epsilon_k,\Sp_n)$};\\
\frac{1}{2}k(n-k-1), & \text{if\/ $(\bfG_k,\bfG'_n)=(\Sp_k,\rmO^\epsilon_n)$};\\
k(n-k), & \text{if\/ $(\bfG_k,\bfG'_n)=(\rmU_k,\rmU_n)$}.
\end{cases}
\]
\end{thm}

\begin{rem}
If $(\bfG_k,\bfG'_n)=(\rmO_k^\epsilon,\Sp_n)$ and the dual pair is in the stable range,
then the $\eta$-correspondence and the $\underline\theta$-correspondence coincide,
and every $\rho\in\cale(G_k)$ is $(n-k)$-admissible (\cf.~Lemma~\ref{0516} and Lemma~\ref{0517}).
For this case, the above formula is already known in \cite{gurevich-howe}, (1.11) and theorem 4.1.1.
\end{rem}

An application of the above theorem is to determine the maximum and minimum values of the set
\[
\{\,\deg_q(\rho)\mid\rho\in\cale(G_n),\ \Theta\text{\rm -rk}(\rho)=k\,\}
\]
where $\bfG_n=\Sp_n$, $\rmO_n^\epsilon$ or $\rmU_n$, and $k\leq n$.
This problem will be investigated in another papers of the author.
A lower bound and an upper bound of the above set can be found in \cite{GLT02} and \cite{GLT01}.

\subsection{}
The contents of this article are as follows.
In Section 2 we set up basic notations and express the formula of $\deg_q(\rho)$
for any unipotent characters in terms of Lusztig symbols.
In Section 3, we investigate the $\Theta$-correspondence and define $\ell$-admissibility
on unipotent characters for a dual pair of a symplectic group and an even orthogonal group.
Moreover, Theorem~\ref{0101} and Theorem~\ref{0501} are proved for this case.
In Section 4, we prove our main results for the cases of $\Theta$-correspondence of
unipotent characters for a dual pair of two unitary groups.
We prove Theorem~\ref{0101} and Theorem~\ref{0501} completely via the Lusztig correspondence
in the final section. 

\section{Degrees in $q$ of Unipotent Characters}

\subsection{Basic notations}
Let $\bfG$ be a classical group defined over $\bfF_q$,
and let $G=\bfG(q)$ denote the finite group of rational points.

For a polynomial $f$ in $q$, let $\deg_q(f)$ denote the degree of the polynomial $f$.
Let $|G|_{p'}$ denote the part of the order $|G|$ prime to $p$.
Now $|G|_{p'}$ is a polynomial in $q$ and it is well-known that
\[
\deg_q(|G|_{p'})=\begin{cases}
\tfrac{1}{4}k(k+2), & \text{if $\bfG=\Sp_k$ or $\rmO_{k+1}$, $k$ even}; \\
\tfrac{1}{4}k^2, & \text{if $\bfG=\rmO^\epsilon_k$, $k$ even}; \\
\tfrac{1}{2}k(k+1), & \text{if $\bfG=\rmU_k$}.
\end{cases}
\]

\subsection{Symbols and unipotent characters}
Let $\mu=[\mu_1,\ldots,\mu_k]$ and $\mu=[\nu_1,\ldots,\nu_l]$ be two partitions.
We may assume that $k=l$ by adding some $0$'s if necessary.
Then we denote
\begin{equation}\label{0306}
\mu\preccurlyeq\nu\quad\text{ if \ }\nu_1\geq\mu_1\geq\nu_2\geq\mu_2\geq\cdots\geq\nu_k\geq\mu_k.
\end{equation}
An ordered pair $\sqbinom{\mu}{\nu}$ of two partitions is called a \emph{bi-partition} of $n$
if
\[
\left|\sqbinom{\mu}{\nu}\right|:=|\mu|+|\nu|=n.
\]
The set of bi-partitions of $n$ is denoted by $\calp_2(n)$.
For a bi-partition $\sqbinom{\mu}{\nu}$, we define its \emph{transpose} by
$\sqbinom{\mu}{\nu}^\rmt=\sqbinom{\nu}{\mu}$.

A $\beta$-set is a finite subset $A=\{a_1,a_2,\ldots,a_m\}$ of non-negative integers
written strictly decreasingly, i.e., $a_1>a_2>\cdots>a_m$.
A symbol $\Lambda=\binom{A}{B}$ is an ordered pair of two $\beta$-sets.
Let $\cals$ denote the set of symbols $\binom{A}{B}$ such that $0\not\in A\cap B$.
For a symbol $\Lambda=\binom{a_1,a_2,\ldots,a_{m_1}}{b_1,b_2,\ldots,b_{m_2}}$,
we define its \emph{defect} and \emph{rank} by
\begin{align}\label{0201}
\begin{split}
{\rm def}(\Lambda) &=m_1-m_2, \\
{\rm rk}(\Lambda) &=\sum_{i=1}^{m_1}a_i+\sum_{j=1}^{m_2}b_j
-\left\lfloor\left(\frac{m_1+m_2-1}{2}\right)^2\right\rfloor.
\end{split}
\end{align}
We also define the \emph{transpose} of a symbol by $\binom{A}{B}^\rmt=\binom{B}{A}$.
There is a mapping $\Upsilon$ from the set of symbols to the set of bi-partitions:
\begin{equation}
\Upsilon\colon\binom{a_1,\ldots,a_{m_1}}{b_1,\ldots,b_{m_2}}\mapsto
\sqbinom{a_1-(m_1-1),a_2-(m_1-2),\ldots,a_{m_1-1}-1,a_{m_1}}
{b_1-(m_2-1),b_2-(m_2-2),\ldots,b_{m_2-1}-1,b_{m_2}}.
\end{equation}
It is easy to check that
\begin{equation}\label{0206}
|\Upsilon(\Lambda)|=\begin{cases}
{\rm rk}(\Lambda)-\tfrac{1}{4}({\rm def}(\Lambda)-1)({\rm def}(\Lambda)+1), & \text{if ${\rm def}(\Lambda)$ is odd};\\
{\rm rk}(\Lambda)-\tfrac{1}{4}({\rm def}(\Lambda))^2, & \text{if ${\rm def}(\Lambda)$ is even}.
\end{cases}
\end{equation}

For a non-negative even integer $k$, we define a set of symbols associated to a symplectic group
or an even orthogonal group:
\begin{align}\label{0205}
\begin{split}
\cals_{\rmO^+_k} &=\{\,\Lambda\in\cals\mid{\rm rk}(\Lambda)=\tfrac{k}{2},\ {\rm def}(\Lambda)\equiv 0\pmod 4\,\};\\
\cals_{\Sp_k} &=\{\,\Lambda\in\cals\mid{\rm rk}(\Lambda)=\tfrac{k}{2},\ {\rm def}(\Lambda)\equiv 1\pmod 4\,\};\\
\cals_{\rmO^-_k} &=\{\,\Lambda\in\cals\mid{\rm rk}(\Lambda)=\tfrac{k}{2},\ {\rm def}(\Lambda)\equiv 2\pmod 4\,\}.
\end{split}
\end{align}
From \cite{lg} it is known that there exists a bijection from $\cals_\bfG$ to the set $\cale(G)_1$
of unipotent characters of $G$ denoted by $\Lambda\mapsto\rho_\Lambda$.

\begin{rem}
The set $\cals_{\Sp_k}$ of symbols considered here is slightly different from the original one
in \cite{lg}.
Please see \cite{pan-uniform} subsection~3.3 to see the equivalence of the two parametrizations.
\end{rem}

Suppose all the entries $a_1,\ldots,a_{m_1},b_1,\ldots,b_{m_2}$ in a symbol
$\Lambda=\binom{a_1,\ldots,a_{m_1}}{b_1,\ldots,b_{m_2}}$ are
$z_1,z_2,\ldots,z_m$ with $z_1\geq z_2\geq\cdots\geq z_m$
where $m=m_1+m_2$.
The following lemma is \cite{pan-eta} lemma 2.12:

\begin{lem}\label{0202}
Let $\Lambda=\binom{a_1,\ldots,a_{m_1}}{b_1,\ldots,b_{m_2}}\in\cals_\bfG$ where
$\bfG=\Sp_k$ or $\rmO^\epsilon_k$ and $k$ is even.
Then $\deg_q(\rho_\Lambda)$ is equal to
\[
\sum_{i=1}^m (m-i)z_i-z_i(z_i+1)
+\begin{cases}
\frac{1}{4}k(k+2)-\frac{1}{24}(m-1)(m-3)(2m-1), & \text{if $m$ is odd};\\
\frac{1}{4}k^2-\frac{1}{24}m(m-2)(2m-5), & \text{if $m$ is even}.
\end{cases}
\]
\end{lem}

\subsection{Unipotent characters of a unitary group}
For a symbol $\Lambda=\binom{a_1,a_2,\ldots,a_{m_1}}{b_1,b_2,\ldots,b_{m_2}}\in\cals$,
we define
\begin{equation}\label{0402}
{\rm rk}_\rmU(\Lambda)
=\sum_{i=1}^{m_1} 2a_i+\sum_{j=1}^{m_2} 2b_j
+\frac{|m_1-m_2|}{2}-\frac{(m_1+m_2)(m_1+m_2-2)}{2}.
\end{equation}
For a non-negative integer $k$,
we let $\cals_{\rmU_k}$ be the set of symbols $\Lambda\in\cals$ such that
\begin{itemize}
\item ${\rm def}(\Lambda)$ is either even and non-negative, or odd and negative;

\item ${\rm rk}_\rmU(\Lambda)=k$.
\end{itemize}
Again, there is a parametrization from $\cals_{\rmU_k}$ to $\cale(\rmU_k(q))_1$
denoted by $\Lambda\mapsto\rho_\Lambda$ (\cf.~\cite{FS} or \cite{pan-eta-unitary}).

For a symbol $\Lambda=\binom{a_1,a_2,\ldots,a_{m_1}}{b_1,b_2,\ldots,b_{m_2}}$,
we define two $\beta$-sets:
\begin{align}\label{0404}
\begin{split}
X^0=X^0_\Lambda &=\begin{cases}
\{2b_1,2b_2,\ldots,2b_{m_2}\}, & \text{if ${\rm def}(\Lambda)$ is even};\\
\{2a_1,2a_2,\ldots,2a_{m_1}\}, & \text{if ${\rm def}(\Lambda)$ is odd},
\end{cases} \\
X^1=X^1_\Lambda &=\begin{cases}
\{2a_1+1,2a_2+1,\ldots,2a_{m_1}+1\}, & \text{if ${\rm def}(\Lambda)$ is even};\\
\{2b_1+1,2b_2+1,\ldots,2b_{m_2}+1\}, & \text{if ${\rm def}(\Lambda)$ is odd}.
\end{cases}
\end{split}
\end{align}
For a $\beta$-set $A$ or two $\beta$-sets $A,B$, we define the following polynomials in $q$:
\begin{align}\label{0204}
\begin{split}
\Delta(A)
&= \prod_{a,a'\in A,\ a>a'}(q^a-q^{a'}), \\
\Theta(A)
&= \prod_{a\in A}\prod_{h=1}^a(q^h-(-1)^h), \\
\Xi(A,B)
&= \prod_{a\in A,\ b\in B}(q^a+q^b).
\end{split}
\end{align}
For $\Lambda=\binom{a_1,a_2,\ldots,a_{m_1}}{b_1,b_2,\ldots,b_{m_2}}\in\cals_{\rmU_k}$,
it is well known that (\cf.~for example, \cite{pan-eta-unitary} proposition 2.15):
\begin{equation}\label{0203}
\rho_\Lambda(1)
= \frac{\Delta(X^0)\Delta(X^1)\Xi(X^0,X^1)|\rmU_k(q)|_{p'}}
{\Theta(X^0)\Theta(X^1)q^{\binom{m_1+m_2-1}{2}+\binom{m_1+m_2-2}{2}+\cdots+\binom{2}{2}}}
\end{equation}
where $\binom{a}{b}$ denotes the binomial coefficients.

\begin{lem}\label{0403}
For $\Lambda\in\cals_{\rmU_k}$,
write $X_\Lambda=X^0_\Lambda\cup X^1_\Lambda=\{z_1,z_2,\ldots,z_m\}$ such that
$z_1>z_2>\cdots>z_m$ where $m=m_1+m_2$.
Then
\[
\deg_q(\rho_\Lambda)
=\sum_{i=1}^m (m-i)z_i-\sum_{i=1}^m \frac{z_i(z_i+1)}{2}+\frac{k(k+1)}{2}-\frac{m(m-1)(m-2)}{6}.
\]
\end{lem}
\begin{proof}
Because now $z_1>z_2>\cdots>z_m$,
from (\ref{0204}) it is not difficult to see that
\begin{align*}
\deg_q(\Delta(X^0)\Delta(X^1)\Pi(X^0,X^1))
&= \sum_{(z_i,z_j)\in X\times X,\ i\neq j}\max(z_i,z_j)
= \sum_{i=1}^m (m-i)z_i, \\
\deg_q(\Theta(X^0)\Theta(X^1))
&= \sum_{i=1}^m \frac{z_i(z_i+1)}{2}.
\end{align*}
Moreover, we have
\begin{align*}
\deg_q(|\rmU_k(q)|_{p'})
&= \frac{k(k+1)}{2}, \\
\binom{m-1}{2}+\binom{m-2}{2}+\cdots+\binom{2}{2}
&= \frac{m(m-1)(m-2)}{6}.
\end{align*}
Then the lemma follows from (\ref{0203}) directly.
\end{proof}

\section{Unipotent Characters of Symplectic or Orthogonal Groups}

\subsection{Finite theta correspondence of unipotent characters}
Let $(\bfG,\bfG')$ be a reductive dual pair of a symplectic group and an even orthogonal group.
For symbols $\Lambda\in\cals_\bfG$ and $\Lambda'\in\cals_{\bfG'}$,
we write $\Upsilon(\Lambda)=\sqbinom{\mu}{\nu}$ and $\Upsilon(\Lambda')=\sqbinom{\mu'}{\nu'}$.
Now we define a relation on the set of symbols:
\begin{align}\label{0313}
\begin{split}
\calb_{\Sp_k,\rmO^+_n} &=\{\,(\Lambda,\Lambda')\in\cals_{\Sp_k}\times\cals_{\rmO^+_n}
\mid\nu\preccurlyeq\mu',\ \nu'\preccurlyeq\mu,
\ {\rm def}(\Lambda')=-{\rm def}(\Lambda)+1\,\};\\
\calb_{\Sp_k,\rmO^-_n} &=\{\,(\Lambda,\Lambda')\in\cals_{\Sp_k}\times\cals_{\rmO^-_n}
\mid\mu\preccurlyeq\nu',\ \mu'\preccurlyeq\nu,
\ {\rm def}(\Lambda')=-{\rm def}(\Lambda)-1\,\}
\end{split}
\end{align}
where both $k,n$ are even.
It is known that the unipotent characters are preserved by the correspondence
$\Theta_{\bfG,\bfG'}$ (\cf.\ \cite{adams-moy} theorem~3.5).
The following result on the $\Theta$-correspondence of unipotent characters
is from \cite{pan-finite-unipotent} theorem 1.8:

\begin{prop}
Let $(\bfG,\bfG')$ be a reductive dual pair of a symplectic group and an even orthogonal group.
Then $(\rho_\Lambda,\rho_{\Lambda})\in\Theta_{\bfG,\bfG'}$ if and only if
$(\Lambda,\Lambda')\in\calb_{\bfG,\bfG'}$.
\end{prop}

\subsection{$\ell$-admissible symbols}\label{0305}
Consider a dual pair $(\bfG,\bfG')=(\Sp_k,\rmO^\epsilon_n)$ or $(\rmO^\epsilon_k,\Sp_n)$
where both $n,k$ are even and $k\leq n$.
Let $\Lambda\in\cals_\bfG$, and write
\begin{equation}\label{0309}
\Lambda=\binom{a_1,a_2,\ldots,a_{m_1}}{b_1,b_2,\ldots,b_{m_2}},\qquad
\Upsilon(\Lambda)=\sqbinom{\mu_1,\mu_2,\ldots,\mu_{m_1}}{\nu_1,\nu_2,\ldots,\nu_{m_2}}.
\end{equation}
Then from (\ref{0205}) we can write
\begin{equation}\label{0314}
{\rm def}(\Lambda)
=\begin{cases}
4d, & \text{if $\bfG=\rmO^+_k$};\\
4d+1, & \text{if $\bfG=\Sp_k$};\\
4d+2, & \text{if $\bfG=\rmO^-_k$}
\end{cases}
\end{equation}
for some $d\in\bbZ$,
and we define an integer
\begin{equation}\label{0312}
\tau=\begin{cases}
\tfrac{n-k}{2}+2d, & \text{if $(\bfG,\bfG')=(\rmO^+_k,\Sp_n)$ or $(\Sp_k,\rmO^+_n)$};\\
\tfrac{n-k}{2}-(2d+1), & \text{if $(\bfG,\bfG')=(\rmO^-_k,\Sp_n)$ or $(\Sp_k,\rmO^-_n)$}.
\end{cases}
\end{equation}
Note that the number $\tau$ depends on $k,n$ and ${\rm def}(\Lambda)$.
If $\tau\geq 0$,
as in \cite{pan-eta},
$\underline\theta(\Lambda)$ is defined to be the symbol in $\cals_{\bfG'}$
such that
\begin{equation}\label{0301}
\Upsilon(\underline\theta(\Lambda))
=\begin{cases}
\sqbinom{\nu_1,\nu_2,\ldots,\nu_{m_2}}{\mu_1,\mu_2,\ldots,\mu_{m_1}}\cup\sqbinom{\tau}{-}, & \text{if $\epsilon=+$};\\
\sqbinom{\nu_1,\nu_2,\ldots,\nu_{m_2}}{\mu_1,\mu_2,\ldots,\mu_{m_1}}\cup\sqbinom{-}{\tau}, & \text{if $\epsilon=-$}.
\end{cases}
\end{equation}
Moreover, if $\tau\geq\nu_1$ when $\epsilon=+$; or $\tau\geq\mu_1$ when $\epsilon=-$, we have
\[
\underline\theta(\Lambda)=
\begin{cases}
\binom{\tau+m_2,b_1,b_2,\ldots,b_{m_2}}{a_1,a_2,\ldots,a_{m_1}}, & \text{if $\epsilon=+$};\\
\binom{b_1,b_2,\ldots,b_{m_2}}{\tau+m_1,a_1,a_2,\ldots,a_{m_1}}, & \text{if $\epsilon=-$}.
\end{cases}
\]

Let $\Lambda\in\cals_\bfG$ with $\Upsilon(\Lambda)$ given as in (\ref{0309})
and ${\rm def}(\Lambda)$ given as in (\ref{0314}),
and let $\ell$ be a non-negative even integer.
Then $\Lambda$ is called \emph{$\ell$-admissible} for the dual pair $(\bfG,\bfG')$ if the following
condition holds:
\begin{equation}\label{0317}
\begin{cases}
\mu_1\leq\frac{\ell}{2}-2d\text{ and }\nu_1\leq\frac{\ell}{2}+2d, & \text{if $(\bfG,\bfG')=(\Sp,\rmO^+_{\rm even})$};\\
\mu_1\leq\frac{\ell}{2}-2d-1\text{ and }\nu_1\leq\frac{\ell}{2}+2d+1, & \text{if $(\bfG,\bfG')=(\Sp,\rmO^-_{\rm even})$};\\
\mu_1\leq\frac{\ell}{2}-2d+1\text{ and }\nu_1\leq\frac{\ell}{2}+2d, & \text{if $(\bfG,\bfG')=(\rmO^+_{\rm even},\Sp)$};\\
\mu_1\leq\frac{\ell}{2}-2d-1\text{ and }\nu_1\leq\frac{\ell}{2}+2d+2, & \text{if $(\bfG,\bfG')=(\rmO^-_{\rm even},\Sp)$}.
\end{cases}
\end{equation}
A unipotent character $\rho_\Lambda\in\cale(G)_1$ is called \emph{$\ell$-admissible}
if $\Lambda$ is $\ell$-admissible.
It is obvious that if $\Lambda$ is $\ell$-admissible, then it is also
$\ell'$-admissible for any even $\ell'\geq \ell$.

\begin{exam}\label{0307}
In this example we consider $0$-admissible symbols $\Lambda\in\cals_\bfG$
for a dual pair $(\bfG,\bfG')$.
Write $\Lambda=\binom{a_1,a_2,\ldots,a_{m_1}}{b_1,b_2,\ldots,b_{m_2}}$ and
$\Upsilon(\Lambda)=\sqbinom{\mu_1,\mu_2,\ldots,\mu_{m_1}}{\nu_1,\nu_2,\ldots,\nu_{m_2}}$.
\begin{enumerate}
\item From (\ref{0317}),
a symbol $\Lambda\in\cals_{\Sp_k}$ of defect $4d+1$ is $0$-admissible for $(\Sp_k,\rmO^+_n)$ if and only if
$\mu_1\leq -2d$ and $\nu_1\leq 2d$.
Because $\mu_1,\nu_1$ are non-negative,
we must have $d=0$ and hence $\mu_1=\nu_1=0$.
Therefore, $\Upsilon(\Lambda)=\sqbinom{0}{0}$, and then $k=0$.

\item A symbol $\Lambda\in\cals_{\Sp_k}$ of defect $4d+1$ is $0$-admissible for $(\Sp_k,\rmO^-_n)$
if and only if $\mu_1\leq -2d-1$ and $\nu_1\leq 2d+1$.
There is no integer $d$ such that both $-2d-1$ and $2d+1$ are non-negative.

\item A symbol $\Lambda\in\cals_{\rmO^+_k}$ of defect $4d$ is $0$-admissible for $(\rmO^+_k,\Sp_n)$
if and only if $\mu_1\leq -2d+1$ and $\nu_1\leq 2d$.
The only possibility is $d=0$ and hence $\mu_1\leq 1$ and $\nu_1\leq 0$, i.e.,
\[
\Upsilon(\Lambda)=\sqbinom{1,1,\ldots,1}{0,0,\ldots,0}\in\calp_2(k),\qquad
\Lambda=\binom{k,k-1,\ldots,1}{k-1,k-2,\ldots,0}.
\]

\item A symbol $\Lambda\in\cals_{\rmO^-_k}$ of defect $4d+2$ is $0$-admissible for $(\rmO^-_k,\Sp_n)$
if and only if $\mu_1\leq -2d-1$ and $\nu_1\leq 2d+2$.
The only possibility is $d=-1$ and hence $\mu_1\leq 1$ and $\nu_1\leq 0$, i.e.,
\[
\Upsilon(\Lambda)=\sqbinom{1,1,\ldots,1}{0,0,\ldots,0}\in\calp_2(k-1),\qquad
\Lambda=\binom{k-1,k-2,\ldots,1}{k,k-1,\ldots,0}.
\]
\end{enumerate}
\end{exam}

\begin{lem}\label{0319}
Consider the dual pair $(\bfG,\bfG')=(\Sp_k,\rmO^\epsilon_n)$ or $(\rmO^\epsilon_k,\Sp_n)$
where both $n,k$ are even and $k\leq n$.
Then every unipotent character $\rho_\Lambda\in\cale(G)_1$ is $k$-admissible.
Consequently, if $(\bfG,\bfG')$ is in stable range,
then every unipotent character $\rho_\Lambda\in\cale(G)_1$ is $(n-k)$-admissible.
\end{lem}
\begin{proof}
Let $\Lambda\in\cals_\bfG$, and write $\Lambda=\binom{a_1,a_2,\ldots,a_{m_1}}{b_1,b_2,\ldots,b_{m_2}}$ 
and $\Upsilon(\Lambda)=\sqbinom{\mu_1,\mu_2,\ldots,\mu_{m_1}}{\nu_1,\nu_2,\ldots,\nu_{m_2}}$.
\begin{enumerate}
\item Suppose that $(\bfG,\bfG')=(\Sp_k,\rmO^+_n)$,
and then ${\rm def}(\Lambda)=4d+1$ for some $d\in\bbZ$.
Then we have
\begin{align*}
\mu_1 &\leq |\Upsilon(\Lambda)|=\frac{k}{2}-2d(2d+1)=\frac{k}{2}-2d-4d^2\leq\frac{k}{2}-2d, \\
\nu_1 &\leq |\Upsilon(\Lambda)|=\frac{k}{2}-2d(2d+1)=\frac{k}{2}+2d-(2d+1)^2\leq\frac{k}{2}+2d.
\end{align*}

\item Suppose that $(\bfG,\bfG')=(\Sp_k,\rmO^-_n)$,
and then ${\rm def}(\Lambda)=4d+1$ for some $d\in\bbZ$.
Then we have
\begin{align*}
\mu_1 &\leq |\Upsilon(\Lambda)|=\frac{k}{2}-2d(2d+1)\leq\frac{k}{2}-2d-1, \\
\nu_1 &\leq |\Upsilon(\Lambda)|=\frac{k}{2}-2d(2d+1)\leq\frac{k}{2}+2d+1.
\end{align*}

\item Suppose that $(\bfG,\bfG')=(\rmO^+_k,\Sp_n)$,
and then ${\rm def}(\Lambda)=4d$ for some $d\in\bbZ$.
Then we have
\begin{align*}
\mu_1 &\leq |\Upsilon(\Lambda)|=\frac{k}{2}-4d^2\leq\frac{k}{2}-2d+1, \\
\nu_1 &\leq |\Upsilon(\Lambda)|=\frac{k}{2}-4d^2\leq\frac{k}{2}+2d.
\end{align*}

\item Suppose that $(\bfG,\bfG')=(\rmO^-_k,\Sp_n)$,
and then ${\rm def}(\Lambda)=4d+2$ for some $d\in\bbZ$.
Then we have
\begin{align*}
\mu_1 &\leq |\Upsilon(\Lambda)|=\frac{k}{2}-(2d+1)^2\leq\frac{k}{2}-2d-1, \\
\nu_1 &\leq |\Upsilon(\Lambda)|=\frac{k}{2}-(2d+1)^2\leq\frac{k}{2}+2d+2.
\end{align*}
\end{enumerate}
For all cases, we see that $\Lambda$ is $k$-admissible.

Now if $(\bfG,\bfG')$ is in stable range, then we have $k\leq n-k$.
Therefore, $\Lambda$ is also $(n-k)$-admissible.
\end{proof}

The following lemma means that if $\Lambda\in\cals_\bfG$ is $(n-k)$-admissible,
then $\underline\theta(\Lambda)$ (\cf.~(\ref{0312})) is defined for the pair
$(\bfG,\bfG')=(\Sp_k,\rmO^\epsilon_n)$ or $(\rmO^\epsilon_k,\Sp_n)$.

\begin{lem}\label{0318}
Consider the dual pair $(\bfG,\bfG')=(\Sp_k,\rmO^\epsilon_n)$ or $(\rmO^\epsilon_k,\Sp_n)$
where both $n,k$ are even and $k\leq n$.
If $\Lambda$ is $(n-k)$-admissible,
then $\tau\geq 0$.
\end{lem}
\begin{proof}
First, suppose that either
\begin{itemize}
\item $(\bfG,\bfG')=(\Sp_k,\rmO^+_n)$ and ${\rm def}(\Lambda)=4d+1$, or

\item $(\bfG,\bfG')=(\rmO^+_k,\Sp_n)$ and ${\rm def}(\Lambda)=4d$
\end{itemize}
for some $d\in\bbZ$.
If $\Lambda$ is $(n-k)$-admissible,
then by (\ref{0312}) and (\ref{0317}) we have $\tau=\frac{1}{2}(n-k)+2d\geq\nu_1\geq 0$.
Next, suppose that either
\begin{itemize}
\item $(\bfG,\bfG')=(\Sp_k,\rmO^-_n)$ and ${\rm def}(\Lambda)=4d+1$, or

\item $(\bfG,\bfG')=(\rmO^-_k,\Sp_n)$ and ${\rm def}(\Lambda)=4d+2$
\end{itemize}
for some $d\in\bbZ$.
If $\Lambda$ is $(n-k)$-admissible,
then we have $\tau=\frac{1}{2}(n-k)-2d-1\geq\mu_1\geq 0$.
\end{proof}

Now we recall the definition of \emph{$\Theta$-rank} of $\rho'\in\cale(G')$,
denoted by $\Theta\text{\rm -rk}(\rho')$, from \cite{pan-theta-rank} subsection~3.1:
\begin{itemize}
\item If $\bfG'$ is a symplectic group, we consider dual pairs $(\rmO^\epsilon_k,\bfG')$ and define
\begin{equation}\label{0310}
\Theta\text{\rm -rk}(\rho')
=\min\{\,k\mid\rho'\in\Theta(\rho)\text{ for some }\rho\in\cale(\rmO^\epsilon_k(q)),\ \epsilon=+\text{ or }-\,\}.
\end{equation}

\item If $\bfG'$ is an orthogonal group, we consider dual pairs $(\Sp_k,\bfG')$ and define
\begin{multline}\label{0311}
\Theta\text{\rm -rk}(\rho')
=\min\{\,k\mid\rho'\chi'\in\Theta(\rho)
\text{ for some }\rho\in\cale(\Sp_k(q)) \\
\text{ and some linear character }\chi'\in\cale(G')\,\}.
\end{multline}
Note that an orthogonal group $G'$ has at most four linear characters:
${\bf 1}$, $\sgn$, $\chi_{\bfG'}$ and $\chi_{\bfG'}\sgn$
where $\chi_{\bfG'}$ denotes the linear character of order two given by the spinor norm.
\end{itemize}

\begin{lem}\label{0316}
Consider the dual pair $(\bfG,\bfG')=(\Sp_k,\rmO^\epsilon_n)$ or $(\rmO^\epsilon_k,\Sp_n)$
where both $n,k$ are even and $k\leq n$.
If $\Lambda\in\cals_\bfG$ is $(n-k)$-admissible,
then $\underline\theta(\rho_\Lambda)$ is of\/ $\Theta$-rank $k$.
\end{lem}
\begin{proof}
Write $\Lambda$ and $\Upsilon(\Lambda)$ as in (\ref{0309}).
Suppose that $\Lambda\in\cals_\bfG$ is $(n-k)$-admissible.
Then $\tau\geq 0$ by Lemma~\ref{0318}, and so the unipotent character
$\underline\theta(\rho_\Lambda)$ is defined.

\begin{enumerate}
\item Suppose that $(\bfG,\bfG')=(\Sp_k,\rmO^+_n)$ and then ${\rm def}(\Lambda)=4d+1$
for some $d\in\bbZ$.
Now the condition $\nu_1\leq \frac{1}{2}(n-k)+2d$ means that $\nu_1\leq\tau$ and then
\[
\Upsilon(\underline\theta(\Lambda))
=\sqbinom{\tau,\nu_1,\nu_2,\ldots,\nu_{m_2}}{\mu_1,\mu_2,\ldots,\mu_{m_1}}.
\]
By the result in \cite{pan-Lusztig-correspondence} section 8,
we see that the unipotent character $\underline\theta(\rho_\Lambda)\in\cale(\rmO^+_n(q))$ does not occur in
the $\Theta$-correspondence for the dual pair $(\Sp_{k'},\rmO^+_n)$ for any $k'<k$.
The condition $\mu_1\leq \frac{n-k}{2}-2d$ means that $\mu_1+(m_1-1)\leq\tau+m_2$
from (\ref{0312}) and (\ref{0301}),
i.e., $a_1\leq \tau+m_2$.
This implies that the unipotent character $\underline\theta(\rho_\Lambda)\sgn\in\cale(\rmO^+_n(q))$
does not occur in the $\Theta$-correspondence for the dual pair $(\Sp_{k''},\rmO^+_n)$ for any $k''<k$.
Moreover, by \cite{pan-theta-rank} Lemma~2.16, we see that both
$\underline\theta(\rho_\Lambda)\chi_{\rmO^+_n},\underline\theta(\rho_\Lambda)\chi_{\rmO^+_n}\sgn$
do not occur in the $\Theta$-correspondence for dual pair $(\Sp_{k'''},\rmO^+_n)$ for any $k'''<k$.

\item Suppose that $(\bfG,\bfG')=(\Sp_k,\rmO^-_n)$ and then ${\rm def}(\Lambda)=4d+1$
for some $d\in\bbZ$.
Then the condition $\mu_1\leq \frac{n-k}{2}-2d-1$ means that $\mu_1\leq\tau$ and then
\[
\Upsilon(\underline\theta(\Lambda))
=\sqbinom{\nu_1,\nu_2,\ldots,\nu_{m_2}}{\tau,\mu_1,\mu_2,\ldots,\mu_{m_1}}.
\]
Hence $\underline\theta(\rho_\Lambda)\in\cale(\rmO^-_n(q))$ does not occur in
the $\Theta$-correspondence for the dual pair $(\Sp_{k'},\rmO^-_n)$ for any $k'<k$.
The condition $\nu_1\leq \frac{n-k}{2}+2d+1$ means that $b_1\leq\tau+m_1$.
This implies that the unipotent character $\underline\theta(\rho_\Lambda)\sgn\in\cale(\rmO^-_n(q))$
does not occur in the $\Theta$-correspondence for the dual pair $(\Sp_{k''},\rmO^-_n)$ for any $k''<k$.
As in (1), we know that both
$\underline\theta(\rho_\Lambda)\chi_{\rmO^-_n},\underline\theta(\rho_\Lambda)\chi_{\rmO^-_n}\sgn$
do not occur in the $\Theta$-correspondence for dual pair $(\Sp_{k'''},\rmO^-_n)$ for any $k'''<k$.

\item Suppose that $(\bfG,\bfG')=(\rmO^+_k,\Sp_n)$ and then ${\rm def}(\Lambda)=4d$
for some $d\in\bbZ$.
Then the conditions $\nu_1\leq \frac{n-k}{2}+2d$ means that $\nu_1\leq\tau$.
Then the unipotent character $\underline\theta(\rho_\Lambda)\in\cale(\Sp_n(q))$ does not occur in
the $\Theta$-correspondence for the dual pair $(\rmO^+_{k'},\Sp_n)$ for any even $k'<k$.
The condition $\mu_1\leq \frac{n-k}{2}-2d+1$ means that $a_1\leq\tau+m_2$.
This implies that the unipotent character $\underline\theta(\rho_\Lambda)\in\cale(\Sp_n(q))$
does not occur in the $\Theta$-correspondence for the dual pair $(\rmO^-_{k''},\Sp_n)$ for any even $k''<k$.
From \cite{pan-theta-rank} lemma~2.13, we know that any unipotent character of $\Sp_n(q)$ does not occur in
the $\Theta$-correspondence for the dual pair $(\rmO_{k'''},\Sp_n)$ for any odd $k'''<n$.

\item Suppose that $(\bfG,\bfG')=(\rmO^-_k,\Sp_n)$ and then ${\rm def}(\Lambda)=4d+2$
for some $d\in\bbZ$.
Then condition $\mu_1\leq \frac{n-k}{2}-2d-1$ means that $\mu_1\leq\tau$,
and the condition $\nu_1\leq \frac{n-k}{2}+2d+2$ means that $b_1\leq\tau+m_2$.
The remaining proof is similar to that of (3).
\end{enumerate}
For all four cases, we conclude that $\underline\theta(\rho_\Lambda)$ is of $\Theta$-rank $k$
from (\ref{0310}) and (\ref{0311}).
\end{proof}

\begin{prop}\label{0302}
Suppose that both $n,k$ are even and $k\leq n$.
Then a unipotent character $\rho'\in\cale(\rmO^\epsilon_n(q))_1$ is of\/ $\Theta$-rank $k$ if and only if
$\rho'$ or $\rho'\sgn$ is equal to $\underline\theta(\rho_\Lambda)$ for some
$(n-k)$-admissible unipotent character $\rho_\Lambda\in\cale(\Sp_k(q))$.
\end{prop}
\begin{proof}
Let $\rho'$ be a unipotent character of $\rmO^\epsilon_n(q)$.
If $\rho'=\underline\theta(\rho_\Lambda)$ or $\rho'\cdot\sgn=\underline\theta(\rho_\Lambda)$
for some $(n-k)$-admissible $\Lambda\in\cals_{\Sp_k}$,
then $\Theta\text{\rm -rk}(\rho')=\Theta\text{\rm -rk}(\rho'\cdot\sgn)=k$ by Lemma~\ref{0316}.

Now suppose that $\Theta\text{\rm -rk}(\rho')=k$.
Then from \cite{pan-eta} and the definition of $\Theta$-rank,
we know that $\rho'$ or $\rho'\cdot\sgn$ is equal to
$\rho_{\underline\theta(\Lambda)}=\underline\theta(\rho_\Lambda)$
for some $\Lambda\in\cals_{\Sp_k}$.
Without loss of generality, we may assume that $\rho'=\underline\theta(\rho_\Lambda)$.
So now we need to show that $\Lambda$ is $(n-k)$-admissible.
Write
\[
\Upsilon(\Lambda)=\sqbinom{\mu_1,\mu_2,\ldots,\mu_{m_1}}{\nu_1,\nu_2,\ldots,\nu_{m_2}}.
\]
Suppose that $\epsilon=+$ and ${\rm def}(\Lambda)=4d+1$ for some $d\in\bbZ$.
Then
\[
\Upsilon(\underline\theta(\Lambda))
=\sqbinom{\nu_1,\nu_2,\ldots,\nu_{m_2}}{\mu_1,\mu_2,\ldots,\mu_{m_1}}\cup\sqbinom{\tau}{-}
\]
by (\ref{0301}).
Note that ${\rm def}(\underline\theta(\Lambda))=-4d$.
If $\Lambda$ is not $(n-k)$-admissible, by definition there are two possibilities:
\begin{itemize}
\item $\nu_1>\frac{1}{2}(n-k)+2d$, or

\item $\mu_1>\frac{1}{2}(n-k)-2d$.
\end{itemize}
Suppose that $\nu_1>\frac{1}{2}(n-k)+2d=\tau$.
Then the rank of the symbol $\Lambda'$ such that
\[
\Upsilon(\Lambda')=\sqbinom{\mu_1,\mu_2,\ldots,\mu_{m_1}}{\nu_2,\nu_3,\ldots,\nu_{m_2}}\cup\sqbinom{-}{\tau}
\]
is less than the rank of $\Lambda$.
Then we see that $(\rho_{\Lambda'},\rho_{\underline\theta(\Lambda)})$ occurs in the correspondence $\Theta$
for a dual pair $(\Sp_{k'},\rmO^+_n)$ for some $k'<k$.

Now $\rho'\cdot\sgn=\rho_{\underline\theta(\Lambda)^\rmt}$ and ${\rm def}(\underline\theta(\Lambda)^\rmt)=4d$.
By the same argument as above (with $d$ replaced by $-d$ and $\nu_1$ replaced by $\mu_1$),
we see that if $\mu_1>\frac{1}{2}(n-k)-2d$, then $\rho'\cdot\sgn$ occurs in the $\Theta$-correspondence
for a dual pair $(\Sp_{k''},\rmO^+_n)$ for some $k''<k$.
Therefore we conclude that $\Theta\text{\rm -rk}(\rho')<k$
if $\Lambda$ is not $(n-k)$-admissible.

The proof for the case that $\epsilon=-$ is similar.
\end{proof}

\begin{prop}\label{0315}
Suppose that both $n,k$ are even and $k\leq n$.
Then a unipotent character $\rho'\in\cale(\Sp_n(q))_1$ is of\/ $\Theta$-rank $k$ if and only if
$\rho'=\underline\theta(\rho_\Lambda)$ for some $(n-k)$-admissible unipotent character
$\rho_\Lambda\in\cale(\rmO^\epsilon_k(q))$ and some $\epsilon=+$ or $-$.
\end{prop}
\begin{proof}
Let $\rho'$ be a unipotent character of $\Sp_n(q)$.
If $\rho'=\underline\theta(\rho_\Lambda)$ for some $(n-k)$-admissible $\Lambda\in\cals_{\rmO^\epsilon_k}$
for some $\epsilon=+$ or $-$,
then $\Theta\text{\rm -rk}(\rho')=k$ by Lemma~\ref{0316}.

Now suppose that $\Theta\text{\rm -rk}(\rho')=k$.
Then we know that $\rho'=\rho_{\underline\theta(\Lambda)}=\underline\theta(\rho_\Lambda)$
for some $\Lambda\in\cals_{\rmO^\epsilon_k}$ and some $\epsilon=+$ or $-$.
Suppose that $\epsilon=+$ and ${\rm def}(\Lambda)=4d$ for some $d\in\bbZ$.
If $\Lambda$ is not $(n-k)$-admissible, then
\begin{itemize}
\item $\mu_1>\frac{1}{2}(n-k)-2d+1$; or

\item $\nu_1>\frac{1}{2}(n-k)+2d$.
\end{itemize}
If $\nu_1>\frac{1}{2}(n-k)+2d=\tau$, then we know that
$\rho_{\underline\theta(\Lambda)}$ occurs in the $\Theta$-correspondence for the dual pair
$(\rmO^+_{k'},\Sp_n)$ for some even $k'<k$.
If $\mu_1>\frac{1}{2}(n-k)-2d+1$, then we know that
$\rho_{\underline\theta(\Lambda)}$ occurs in the $\Theta$-correspondence for the dual pair
$(\rmO^-_{k''},\Sp_n)$ for some even $k''<k$.
Therefore $\Lambda$ must be $(n-k)$-admissible.
The proof for case $\epsilon=-$ is similar.
\end{proof}

\subsection{Degree difference for symplectic/orthogonal dual pairs}

\begin{prop}\label{0303}
Consider the dual pair $(\bfG,\bfG')=(\rmO^\epsilon_k,\Sp_n)$ or $(\Sp_k,\rmO^\epsilon_n)$
where both $k,n$ are even and $k\leq n$, $\epsilon=+$ or $-$.
If $\Lambda\in\cals_\bfG$ is $(n-k)$-admissible,
then
\[
\deg_q(\underline\theta(\rho_\Lambda))
= \deg_q(\rho_\Lambda)+
\begin{cases}
\frac{1}{2}k(n-k+1), & \text{if\/ $(\bfG,\bfG')=(\rmO^\epsilon_k,\Sp_n)$};\\
\frac{1}{2}k(n-k-1), & \text{if\/ $(\bfG,\bfG')=(\Sp_k,\rmO^\epsilon_n)$}.
\end{cases}
\]
\end{prop}
\begin{proof}
Let $\Lambda\in\cals_\bfG$ be $(n-k)$-admissible,
and write $\Lambda=\binom{a_1,a_2,\ldots,a_{m_1}}{b_1,b_2,\ldots,b_{m_2}}$,
$\Upsilon(\Lambda)=\sqbinom{\mu_1,\mu_2,\ldots,\mu_{m_1}}{\nu_1,\nu_2,\ldots,\nu_{m_2}}$.
Let $z_1,z_2,\ldots,z_m$ denote the entries $a_1,\ldots,a_{m_1},b_1,\ldots,b_{m_2}$ of $\Lambda$
where $z_1\geq z_2\geq\cdots\geq z_m$ and $m=m_1+m_2$.
By Lemma~\ref{0318} and Lemma~\ref{0316}, we know that $\underline\theta(\rho_\Lambda)$ is defined and
of $\Theta$-rank $k$.
Let $\{z'_1,z'_2,\ldots,z'_{m+1}\}$ denote the set of entries of $\underline\theta(\Lambda)$
where $z'_1\geq z'_2\geq\cdots\geq z'_{m+1}$.
From (\ref{0201}), we know that
\begin{align}\label{0304}
\begin{split}
\frac{k}{2}={\rm rk}(\Lambda)
&= \sum_{i=1}^m z_i-\left\lfloor\left(\frac{m-1}{2}\right)^2\right\rfloor, \\
\frac{n}{2}={\rm rk}(\underline\theta(\Lambda))
&= \sum_{i=1}^{m+1} z'_i-\left\lfloor\left(\frac{m}{2}\right)^2\right\rfloor.
\end{split}
\end{align}
Now we consider the following cases:
\begin{enumerate}
\item Suppose that $(\bfG,\bfG')=(\rmO^+_k,\Sp_n)$.
Now by (\ref{0313}) and (\ref{0301}) we can write
$\underline\theta(\Lambda)=\binom{b'_1,b'_2,\ldots,b'_{m_2+1}}{a_1,a_2,\ldots,a_{m_1}}$ for some $b'_i$.
The assumption that $\Lambda$ is $(n-k)$-admissible
implies that $b'_1=\tau+m_2$, $b'_{j+1}=b_j$ for $j=1,2,\ldots,m_2$, and $b'_1\geq a_1$.
Hence $z_1'=b_1'$ and $z'_{i+1}=z_i$ for $i=1,2,\ldots,m$.
Now ${\rm def}(\Lambda)\equiv 0\pmod 4$, so $m$ is even.
From (\ref{0304}), we obtain
\[
z'_1
=\frac{n}{2}+\left(\frac{m}{2}\right)^2-\frac{k}{2}-\left\lfloor\left(\frac{m-1}{2}\right)^2\right\rfloor
=\frac{1}{2}(n-k+m).
\]
Now $\deg_q(|\rmO^\epsilon_k(q)|_{p'})=\frac{1}{4}k^2$ and
$\deg_q(|\Sp_n(q)|_{p'})=\frac{1}{4}n(n+2)$.
Moreover,
\[
(m+1-(i+1))z'_{i+1}-z'_{i+1}(z'_{i+1}+1)
=(m-i)z_i-z_i(z_i+1)
\]
for $i=1,2,\ldots,m$.
Then by Lemma~\ref{0202}, we have
\begin{align*}
& \deg_q(\underline\theta(\rho_{\Lambda}))-\deg_q(\rho_\Lambda) \\
&= m z'_1-z'_1(z'_1+1)-\tfrac{1}{4}m(m-2)+\tfrac{1}{4}n(n+2)-\tfrac{1}{4}k^2 \\
&= \tfrac{1}{2}m(n-k+m)-\tfrac{1}{4}(n-k+m)(n-k+m+2)-\tfrac{1}{4}m(m-2)+\tfrac{1}{4}n(n+2)-\tfrac{1}{4}k^2 \\
&= \tfrac{1}{2}k(n-k+1).
\end{align*}

\item Suppose that $(\bfG,\bfG')=(\rmO^-_k,\Sp_n)$.
Now we can write
$\underline\theta(\Lambda)=\binom{b_1,b_2,\ldots,b_{m_2}}{a'_1,a'_2,\ldots,a'_{m_1+1}}$.
The assumption that $\Lambda$ is $(n-k)$-admissible
implies that $a'_1=\tau+m_1$, $a'_{j+1}=a_j$ for $j=1,2,\ldots,m_1$, and $a'_1\geq b_1$.
Hence $z_1'=a_1'$ and $z'_{i+1}=z_i$ for $i=1,2,\ldots,m$.
Then the remaining argument is exactly the same as in the above case.

\item Suppose that $(\bfG,\bfG')=(\Sp_k,\rmO^\epsilon_n)$.
Now $\underline\theta(\Lambda)$ is written as in case (1) if $\epsilon=+$;
and as in case (2) if $\epsilon=-$.
Now ${\rm def}(\Lambda)\equiv 1\pmod 4$, so $m$ is odd.
The assumption that $\Lambda$ is $(n-k)$-admissible implies that
\[
z'_1
=\frac{n}{2}+\left\lfloor\left(\frac{m}{2}\right)^2\right\rfloor-\frac{k}{2}-\left(\frac{m-1}{2}\right)^2
=\frac{1}{2}(n-k+m-1)
\]
and $z'_{i+1}=z_i$ for $i=1,2,\ldots,m$.
Then by Lemma~\ref{0202}, we have
\begin{align*}
\deg_q(\underline\theta(\rho_\Lambda))-\deg_q(\rho_\Lambda)
&= m z'_1-z'_1(z'_1+1)-\tfrac{1}{4}(m-1)^2+\tfrac{1}{4}n^2-\tfrac{1}{4}k(k+2) \\
&= \tfrac{1}{2}m(n-k+m-1)-\tfrac{1}{4}(n-k+m-1)(n-k+m+1) \\
&\qquad\qquad\qquad\qquad\qquad\qquad -\tfrac{1}{4}(m-1)^2+\tfrac{1}{4}n^2-\tfrac{1}{4}k(k+2) \\
&= \tfrac{1}{2}k(n-k-1).
\end{align*}
\end{enumerate}
\end{proof}

\section{Unipotent Characters of Unitary Groups}

\subsection{Finite theta correspondence of unipotent characters}

Let $(\bfG,\bfG')=(\rmU_k,\rmU_n)$ be a reductive dual pair of two unitary groups.
For symbols $\Lambda\in\cals_\bfG$ and $\Lambda'\in\cals_{\bfG'}$,
we write $\Upsilon(\Lambda)=\sqbinom{\mu}{\nu}$ and $\Upsilon(\Lambda')=\sqbinom{\mu'}{\nu'}$.
Now we define $\calb_{\rmU_k,\rmU_n}$ as follows:
\begin{itemize}
\item if $k+n$ is even, then
\begin{multline*}
\calb_{\rmU_k,\rmU_n}
=\Biggl\{\,(\Lambda,\Lambda')\in\cals_{\rmU_k}\times\cals_{\rmU_n}
\mid\nu\preccurlyeq\mu',\ \nu'\preccurlyeq\mu,\\
{\rm def}(\Lambda')=\begin{cases}
0, & \text{if ${\rm def}(\Lambda)=0$};\\
-{\rm def}(\Lambda)+1, & \text{if ${\rm def}(\Lambda)\neq 0$}
\end{cases}\,\Biggr\}.
\end{multline*}
\item if $k+n$ is odd, then
\[
\calb_{\rmU_k,\rmU_n}
=\{\,(\Lambda,\Lambda')\in\cals_{\rmU_k}\times\cals_{\rmU_n}
\mid\mu\preccurlyeq\nu',\ \mu'\preccurlyeq\nu,
\ {\rm def}(\Lambda')=-{\rm def}(\Lambda)-1\,\}.
\]
\end{itemize}
Here the relation $\mu\preccurlyeq\nu'$ is given in (\ref{0306}).
It is known that the unipotent characters are preserved by the $\Theta_{\bfG,\bfG'}$
(\cf.\ \cite{adams-moy} theorem~3.5).
The following proposition on the $\Theta$-correspondence of unipotent characters
for a unitary dual pair
is rephrased from \cite{amr} th\'eor\`eme 5.15 (see also \cite{pan-Lusztig-correspondence} proposition 5.13):

\begin{prop}
Let $\rho_\Lambda\in\cale(\rmU_k(q))_1$ and $\rho_{\Lambda'}\in\cale(\rmU_n(q))_1$.
Then $(\rho_\Lambda,\rho_{\Lambda'})\in\Theta_{\rmU_k,\rmU_n}$ if and only if
$(\Lambda,\Lambda')\in\calb_{\rmU_k,\rmU_n}$.
\end{prop}

\subsection{$\ell$-admissible unipotent characters}\label{0415}
Consider the dual pair $(\bfG,\bfG')=(\rmU_k,\rmU_n)$ such that $k\leq n$.
Let $\Lambda\in\cals_{\rmU_k}$, and write
\[
\Lambda=\binom{a_1,a_2,\ldots,a_{m_1}}{b_1,b_2,\ldots,b_{m_2}},\qquad
\Upsilon(\Lambda)=\sqbinom{\mu_1,\mu_2,\ldots,\mu_{m_1}}{\nu_1,\nu_2,\ldots,\nu_{m_2}}.
\]
Suppose that $(\Lambda,\Lambda')\in\calb_{\rmU_k,\rmU_n}$.
Let $d=|{\rm def}(\Lambda)|$ and $d'=|{\rm def}(\Lambda')|$.
Define
\begin{equation}\label{0406}
\tau=\frac{1}{2}\left[\left(n-\frac{d'(d'+1)}{2}\right)-\left(k-\frac{d(d+1)}{2}\right)\right].
\end{equation}
Then from the definition of $\calb_{\rmU_k,\rmU_n}$ we can check that
\begin{equation}
\tau=\begin{cases}
\frac{1}{2}(n-k+d), & \text{if $k+n+d$ is even};\\
\frac{1}{2}(n-k-d-1), & \text{if $k+n+d$ is odd}.
\end{cases}
\end{equation}
Note that $\tau$ is an integer depending on $k,n$ and ${\rm def}(\Lambda)$.
If $\tau\geq 0$,
then $\underline\theta(\Lambda)$ is defined to be the symbol in $\cals_{\rmU_n}$ such that
\[
\Upsilon(\underline\theta(\Lambda))
=\begin{cases}
\Upsilon(\Lambda)^\rmt\cup\sqbinom{\tau}{-}, & \text{if $k+n$ is even};\\
\Upsilon(\Lambda)^\rmt\cup\sqbinom{-}{\tau}, & \text{if $k+n$ is odd}.
\end{cases}
\]

For a non-negative integer $\ell$,
a symbol $\Lambda\in\cals_{\rmU_k}$ is called \emph{$\ell$-admissible} if the following conditions hold:
\begin{equation}\label{0409}
\begin{cases}
\mu_1\leq\frac{1}{2}(\ell-d)\text{ and }\nu_1\leq\frac{1}{2}(\ell+d),
& \text{if $\ell$ is even, $d$ is even};\\
\mu_1\leq\frac{1}{2}(\ell+d+1)\text{ and }\nu_1\leq\frac{1}{2}(\ell-d-1),
& \text{if $\ell$ is even, $d$ is odd};\\
\mu_1\leq\frac{1}{2}(\ell-d-1)\text{ and }\nu_1\leq\frac{1}{2}(\ell+d+1),
& \text{if $\ell$ is odd, $d$ is even};\\
\mu_1\leq\frac{1}{2}(\ell+d)\text{ and }\nu_1\leq\frac{1}{2}(\ell-d),
& \text{if $\ell$ is odd, $d$ is odd}.
\end{cases}
\end{equation}
It is clear that if $\Lambda$ is $\ell$-admissible, then it is also $\ell'$-admissible
for any $\ell'\geq\ell$.
A unipotent character $\rho_\Lambda\in\cale(\rmU_k(q))_1$ is called \emph{$\ell$-admissible}
if $\Lambda\in\cals_{\rmU_k}$ is $\ell$-admissible.

\begin{lem}\label{0407}
Every unipotent character $\rho_\Lambda\in\cale(\rmU_k(q))_1$ is $k$-admissible.
Consequently, if the dual pair $(\bfG,\bfG')=(\rmU_k,\rmU_n)$ is in stable range,
then every $\rho_\Lambda\in\cale(G)_1$ is $(n-k)$-admissible.
\end{lem}
\begin{proof}
Let $\Lambda\in\cals_{\rmU_k}$, $d=|{\rm def}(\Lambda)|$, 
and write $\Lambda=\binom{a_1,a_2,\ldots,a_{m_1}}{b_1,b_2,\ldots,b_{m_2}}$,
$\Upsilon(\Lambda)=\sqbinom{\mu_1,\mu_2,\ldots,\mu_{m_1}}{\nu_1,\nu_2,\ldots,\nu_{m_2}}$.
\begin{enumerate}
\item Suppose that $k$ is even.
Then
\begin{align*}
\mu_1 &\leq |\Upsilon(\Lambda)|=\frac{1}{2}(k-d(d+1))\leq\begin{cases}
\frac{1}{2}(k-d), & \text{if $d$ is even};\\
\frac{1}{2}(k+d+1), & \text{if $d$ is odd},
\end{cases}  \\
\nu_1 &\leq |\Upsilon(\Lambda)|=\frac{1}{2}(k-d(d+1))\leq\begin{cases}
\frac{1}{2}(k+d), & \text{if $d$ is even};\\
\frac{1}{2}(k-d-1), & \text{if $d$ is odd}.
\end{cases}
\end{align*}

\item Suppose that $k$ is odd.
Then
\begin{align*}
\mu_1 &\leq |\Upsilon(\Lambda)|=\frac{1}{2}(k-d(d+1))\leq\begin{cases}
\frac{1}{2}(k-d-1), & \text{if $d$ is even};\\
\frac{2}{2}(k+d), & \text{if $d$ is odd},
\end{cases}  \\
\nu_1 &\leq |\Upsilon(\Lambda)|=\frac{1}{2}(k-d(d+1))\leq\begin{cases}
\frac{1}{2}(k+d+1), & \text{if $d$ is even};\\
\frac{2}{2}(k-d), & \text{if $d$ is odd}.
\end{cases}
\end{align*}
\end{enumerate}
For two cases, $\Lambda$ is always $k$-admissible.

Now if $(\rmU_k,\rmU_n)$ is in stable range, then we have $k\leq n-k$.
Therefore, $\Lambda$ is also $(n-k)$-admissible.
\end{proof}

\begin{lem}\label{0410}
For a non-negative integer $d$,
let $\rho_d=\rho_{\Lambda_d}$ denote the unipotent cuspidal character of\/ $\rmU_k(q)$ where
$k=\frac{1}{2}d(d+1)$.
Then $\Lambda_d$ is $\ell$-admissible for $\ell\geq d$ if $\ell-d$ is even;
and is $\ell$-admissible for $\ell\geq d+1$ if $\ell-d$ is odd.
\end{lem}
\begin{proof}
It is known that $\Lambda_d=\binom{d-1,d-2,\ldots,0}{-}$ if $d$ is even; and
$\Lambda_d=\binom{-}{d-1,d-2,\ldots,0}$ if $d$ is odd,
and $\Upsilon(\Lambda_d)=\sqbinom{0}{0}$, i.e., $\mu_1=\nu_1=0$.
Note that ${\rm def}(\Lambda_d)=d$ if $d$ is even; and ${\rm def}(\Lambda_d)=-d$ if $d$ is odd.
Then the lemma follows from (\ref{0409}) immediately.
\end{proof}

\begin{exam}
Let ${\rm St}_{\rmU_k}$ denote the Steinberg character of $\rmU_k(q)$ where $k\geq 1$.
It is known that $\deg_q({\rm St}_{\rmU_k})=\frac{1}{2}k(k-1)$.
Moreover, the symbol $\Lambda$ associated to ${\rm St}_{\rmU_k}$ satisfies
\[
\Lambda=\begin{cases}
\binom{\frac{k-2}{2},\frac{k-4}{2},\ldots,0}{\frac{k}{2},\frac{k-2}{2},\ldots,1} \\
\binom{\frac{k-1}{2},\frac{k-3}{2},\ldots,1}{\frac{k-1}{2},\frac{k-3}{2},\ldots,0}
\end{cases}\
\Upsilon(\Lambda)=\begin{cases}
\sqbinom{-}{1,1,\ldots,1}, & \text{if $k$ is even};\\
\sqbinom{1,1,\ldots,1}{-}, & \text{if $k$ is odd}.
\end{cases}
\]
Note that now ${\rm def}(\Lambda)=0$ if $k$ is even, ${\rm def}(\Lambda)=-1$ if $k$ is odd.
Therefore we see that ${\rm St}_{\rmU_k}$ is $1$-admissible,
and hence it is $\ell$-admissible for any $\ell\geq 1$.
\end{exam}

\begin{lem}\label{0414}
Consider the dual pair $(\rmU_k,\rmU_n)$ where $k\leq n$.
If $\Lambda\in\cals_{\rmU_k}$ is $(n-k)$-admissible,
then $\tau\geq 0$.
\end{lem}
\begin{proof}
Write $\Upsilon(\Lambda)=\sqbinom{\mu_1,\mu_2,\ldots,\mu_{m_1}}{\nu_1,\nu_2,\ldots,\nu_{m_2}}$.
Suppose that $\Lambda\in\cals_{\rmU_k}$ is $(n-k)$-admissible.
First suppose that $k+n+d$ is even.
The condition that $\Lambda$ is $(n-k)$-admissible means that $(n-k)+d\geq 0$ from (\ref{0409}).
This is exactly the condition that $\tau\geq 0$.

Next, suppose that $k+n+d$ is odd.
The condition that $\Lambda$ is $(n-k)$-admissible means that $(n-k)-d-1\geq 0$ from (\ref{0409}).
Then again we have $\tau\geq 0$.
\end{proof}

Recall that from \cite{pan-theta-rank} subsection~3.1,
the \emph{$\Theta$-rank} of $\rho'\in\cale(\rmU_n(q))$ is defined by
\begin{multline}\label{0405}
\Theta\text{\rm -rk}(\rho')
=\min\{\,k\mid\rho'\chi'\in\Theta(\rho)\text{ for some }\rho\in\cale(\rmU_k(q))\\
\text{and some linear character }\chi'\in\cale(\rmU_n(q))\,\}.
\end{multline}

\begin{lem}\label{0413}
Consider the dual pair $(\bfG,\bfG')=(\rmU_k,\rmU_n)$ where $k\leq n$.
If\/ $\Lambda\in\cals_{\rmU_k}$ is $(n-k)$-admissible,
then $\underline\theta(\rho_\Lambda)$ is of\/ $\Theta$-rank $k$.
\end{lem}
\begin{proof}
Suppose that $\Lambda\in\cals_{\rmU_k}$ is $(n-k)$-admissible for $(\rmU_k,\rmU_n)$.
From the previous lemma, we know that $\tau\geq 0$ and so $\underline\theta(\rho_\Lambda)$ is defined.
Let $d=|{\rm def}(\Lambda)|$.
\begin{enumerate}
\item Suppose that $k+n$ is even.
Now we have $\tau\geq\nu_1$.
Then
\[
\underline\theta(\Lambda)=\binom{\tau+m_2,b_1,b_2,\ldots,b_{m_2}}{a_1,a_2,\ldots,a_{m_1}},\qquad
\Upsilon(\underline\theta(\Lambda))=\sqbinom{\tau,\nu_1,\nu_2,\ldots,\nu_{m_2}}{\mu_1,\mu_2,\ldots,\mu_{m_1}}.
\]
\begin{enumerate}
\item Suppose that $d$ is even.
Then the condition $\nu_1\leq\frac{1}{2}(n-k+d)$ implies that $\rho_{\underline\theta(\Lambda)}$
does not occur in the $\Theta$-correspondence for the dual pair
$(\rmU_{k'},\rmU_n)$ for any $k'<k$ such that $k'+n$ is even.
The condition $\mu_1\leq\frac{1}{2}(n-k-d)$ means that $2(\mu_1+m_1-1)\leq 2(\tau+m_2)-1$,
i.e.,
\[
\text{\rm rk}_\rmU\binom{a_1,a_2,\ldots,a_{m_1}}{b_1,b_2,\ldots,b_{m_2}}
\leq\text{\rm rk}_\rmU\binom{a_2,\ldots,a_{m_1}}{\tau+m_2,b_1,b_2,\ldots,b_{m_2}}.
\]
Therefore the unipotent character $\underline\theta(\rho_\Lambda)$ of $\rmU_n(q)$
does not occur in the $\Theta$-correspondence
for the dual pair $(\rmU_{k''},\rmU_n)$ for any $k''<k$ such that $k''+n$ is odd.

\item Suppose that $d$ is odd.
The condition $\nu_1\leq\frac{1}{2}(n-k-d-1)$ implies that $\underline\theta(\rho_\Lambda)$
does not occur in the $\Theta$-correspondence for the dual pair
$(\rmU_{k'},\rmU_n)$ for any $k'<k$ such that $k'+n$ is even.
The condition $\mu_1\leq\frac{1}{2}(n-k+d+1)$ implies $\underline\theta(\rho_\Lambda)$
does not occur in the $\Theta$-correspondence
for the dual pair $(\rmU_{k''},\rmU_n)$ for any $k''<k$ such that $k''+n$ is odd.
\end{enumerate}

\item Suppose that $k+n$ is odd.
Now we have $\tau\geq\mu_1$.
Then
\[
\underline\theta(\Lambda)=\binom{b_1,b_2,\ldots,b_{m_2}}{\tau+m_1,a_1,a_2,\ldots,a_{m_1}},\qquad
\Upsilon(\underline\theta(\Lambda))=\sqbinom{\nu_1,\nu_2,\ldots,\nu_{m_2}}{\tau,\mu_1,\mu_2,\ldots,\mu_{m_1}}.
\]

\begin{enumerate}
\item Suppose that $d$ is even.
The condition $\mu_1\leq\frac{1}{2}(n-k-d-1)$ implies that $\underline\theta(\rho_\Lambda)$
does not occur in the $\Theta$-correspondence for the dual pair
$(\rmU_{k'},\rmU_n)$ for any $k'<k$ such that $k'+n$ is even.
The condition $\nu_1\leq\frac{1}{2}(n-k+d+1)$ means $2(\nu_1+m_2-1)\leq 2(\tau+m_1)+1$,
i.e., .
\[
\text{\rm rk}_\rmU\binom{a_1,a_2,\ldots,a_{m_1}}{b_1,b_2,\ldots,b_{m_2}}
\leq\text{\rm rk}_\rmU\binom{\tau+m_1,a_1,a_2,\ldots,a_{m_1}}{b_2,\ldots,b_{m_2}}.
\]
Therefore
$\underline\theta(\rho_\Lambda)$ does not occur in the $\Theta$-correspondence
for the dual pair $(\rmU_{k''},\rmU_n)$ for any $k''<k$ such that $k''+n$ is odd.

\item Suppose that $d$ is odd.
The condition $\mu_1\leq\frac{1}{2}(n-k+d)$ implies that $\underline\theta(\rho_\Lambda)$
does not occur in the $\Theta$-correspondence for the dual pair
$(\rmU_{k'},\rmU_n)$ for any $k'<k$ such that $k'+n$ is even.
The condition $\nu_1\leq\frac{1}{2}(n-k-d)$ implies that $\underline\theta(\rho_\Lambda)$ does not occur
in the $\Theta$-correspondence for the dual pair $(\rmU_{k''},\rmU_n)$ for any $k''<k$ such that $k''+n$ is odd.
\end{enumerate}
\end{enumerate}
For all cases, we conclude that $\underline\theta(\rho_\Lambda)$ is of $\Theta$-rank $k$.
\end{proof}

\begin{prop}
Consider the dual pair $(\bfG,\bfG')=(\rmU_k,\rmU_n)$ where $k\leq n$.
A unipotent character $\rho'\in\cale(\rmU_n(q))_1$ is of\/ $\Theta$-rank $k$ if and only if
$\rho'=\underline\theta(\rho_\Lambda)$ for some $(n-k)$-admissible symbol $\Lambda\in\cals_{\rmU_k}$.
\end{prop}
\begin{proof}
Let $\rho'$ be a unipotent character of $\rmU_n(q)$.
Suppose that $\rho'=\underline\theta(\rho_\Lambda)$ for some $(n-k)$-admissible $\Lambda\in\cals_{\rmU_k}$.
Then $\Theta\text{\rm -rk}(\rho')=k$ by Lemma~\ref{0413}.

Now we suppose that $\rho'$ is of $\Theta$-rank $k$.
From \cite{pan-theta-rank}, we know that $\rho'=\underline\theta(\rho_\Lambda)$
for some symbol $\Lambda\in\cals_{\rmU_k}$.
We need to show that $\Lambda$ is $(n-k)$-admissible.
Write $\Upsilon(\Lambda)=\sqbinom{\mu_1,\mu_2,\ldots,\mu_{m_1}}{\nu_1,\nu_2,\ldots,\nu_{m_2}}$
and let $d=|{\rm def}(\Lambda)|$.
From \cite{pan-Lusztig-correspondence} section 8, we can conclude the following.
\begin{enumerate}
\item Suppose that $k+n$ is even and $d$ is even.
\begin{enumerate}
\item If $\mu_1>\frac{1}{2}(n-k-d)$, we know that $\underline\theta(\rho_\Lambda)$ occurs
in the $\Theta$-correspondence for the dual pair $(\rmU_{k'},\rmU_n)$ for some $k'<k$
such that $k'+n$ is even.

\item If $\nu_1>\frac{1}{2}(n-k+d)$, then $\underline\theta(\rho_\Lambda)$ occurs
in the $\Theta$-correspondence for the dual pair $(\rmU_{k''},\rmU_n)$ for some $k''<k$
such that $k''+n$ is odd.
\end{enumerate}

\item Suppose that $k+n$ is even and $d$ is odd.
\begin{enumerate}
\item If $\mu_1>\frac{1}{2}(n-k+d+1)$, we know that $\underline\theta(\rho_\Lambda)$ occurs
in the $\Theta$-correspondence for the dual pair $(\rmU_{k'},\rmU_n)$ for some $k'<k$
such that $k'+n$ is even.

\item If $\nu_1>\frac{1}{2}(n-k-d-1)$, then $\underline\theta(\rho_\Lambda)$ occurs
in the $\Theta$-correspondence for the dual pair $(\rmU_{k''},\rmU_n)$ for some $k''<k$
such that $k''+n$ is odd.
\end{enumerate}

\item Suppose that $k+n$ is odd and $d$ is even.
\begin{enumerate}
\item If $\mu_1>\frac{1}{2}(n-k-d-1)$, we know that $\underline\theta(\rho_\Lambda)$ occurs
in the $\Theta$-correspondence for the dual pair $(\rmU_{k'},\rmU_n)$ for some $k'<k$
such that $k'+n$ is even.

\item If $\nu_1>\frac{1}{2}(n-k+d+1)$, then $\underline\theta(\rho_\Lambda)$ occurs
in the $\Theta$-correspondence for the dual pair $(\rmU_{k''},\rmU_n)$ for some $k''<k$
such that $k''+n$ is odd.
\end{enumerate}

\item Suppose that $k+n$ is odd and $d$ is odd.
\begin{enumerate}
\item If $\mu_1>\frac{1}{2}(n-k+d)$, we know that $\underline\theta(\rho_\Lambda)$ occurs
in the $\Theta$-correspondence for the dual pair $(\rmU_{k'},\rmU_n)$ for some $k'<k$
such that $k'+n$ is even.

\item If $\nu_1>\frac{1}{2}(n-k-d)$, then $\underline\theta(\rho_\Lambda)$ occurs
in the $\Theta$-correspondence for the dual pair $(\rmU_{k''},\rmU_n)$ for some $k''<k$
such that $k''+n$ is odd.
\end{enumerate}
\end{enumerate}
Therefore, for all cases, if $\Lambda\in\cals_{\rmU_k}$ is not $(n-k)$-admissible,
we will conclude that $\Theta\text{\rm -rk}(\underline\theta(\rho_\Lambda))<k$.
\end{proof}

\subsection{Degree differences for unitary dual pairs}
\begin{prop}\label{0401}
Consider the dual pair $(\bfG,\bfG')=(\rmU_k,\rmU_n)$ where $k\leq n$.
If $\Lambda\in\cals_{\rmU_k}$ is $(n-k)$-admissible,
then
\[
\deg_q(\underline\theta(\rho_\Lambda))=\deg_q(\rho_\Lambda)+k(n-k).
\]
\end{prop}
\begin{proof}
Let $\Lambda\in\cals_{\rmU_k}$ be $(n-k)$-admissible, and write
$\Lambda=\binom{a_1,a_2,\ldots,a_{m_1}}{b_1,b_2,\ldots,b_{m_2}}$,
$\Upsilon(\Lambda)=\sqbinom{\mu_1,\mu_2,\ldots,\mu_{m_1}}{\nu_1,\nu_2,\ldots,\nu_{m_2}}$.
\begin{enumerate}
\item Suppose that $k+n$ is odd.
The assumption that $\Lambda$ is $(n-k)$-admissible means that
\[
\underline\theta(\Lambda)=\binom{b_1,b_2,\ldots,b_{m_2}}{\tau+m_1,a_1,a_2,\ldots,a_{m_1}},\qquad
\Upsilon(\underline\theta(\Lambda))=\sqbinom{\nu_1,\nu_2,\ldots,\nu_{m_2}}{\tau,\mu_1,\mu_2,\ldots,\mu_{m_1}}.
\]

\begin{enumerate}
\item Suppose that $m_1+m_2$ is even.
Then ${\rm def}(\Lambda)=m_1-m_2$ is even and non-negative, and
${\rm def}(\underline\theta(\Lambda))=m_2-m_1-1$ is odd and negative.
By (\ref{0402}), we have
\begin{align*}
k &=\sum_{i=1}^{m_1}2a_i+\sum_{j=1}^{m_2}2b_j-\frac{1}{2}m(m-2)+\frac{1}{2}(m_1-m_2) \\
n &=2(\tau+m_1)+\sum_{i=1}^{m_1}2a_i+\sum_{j=1}^{m_2}2b_j-\frac{1}{2}(m+1)(m-1)+\frac{1}{2}(m_1-m_2+1)
\end{align*}
where $m=m_1+m_2$.
Then we have
\[
2(\tau+m_1)+1=n-k+m.
\]
Now by (\ref{0404}), we have
\begin{align*}
X^0_\Lambda&=\{2b_1,2b_2,\ldots,2b_{m_2}\},
& X^1_\Lambda&=\{2a_1+1,2a_2+1,\ldots,2a_{m_1}+1\}; \\
X^0_{\underline\theta(\Lambda)}&=\{2b_1,2b_2,\ldots,2b_{m_2}\},
& X^1_{\underline\theta(\Lambda)}&=\{2(\tau+m_1)+1,2a_1+1,2a_2+1,\ldots,2a_{m_1}+1\}.
\end{align*}
Write
\begin{align*}
X^0_\Lambda\cup X^1_\Lambda =\{z_1,z_2,\ldots,z_m\},\qquad
X^0_{\underline\theta(\Lambda)}\cup X^1_{\underline\theta(\Lambda)} =\{z'_1,z'_2,\ldots,z'_{m+1}\}
\end{align*}
with $z_1>z_2>\cdots>z_m$ and $z'_1>z'_2>\cdots>z'_{m+1}$.
Now $z'_{i+1}=z_i$ for $i=1,\ldots,m$, and $z'_1=n-k+m$.
Hence by Lemma~\ref{0403},
we have
\begin{align*}
&\deg_q(\underline\theta(\rho_\Lambda))-\deg_q(\rho_\Lambda) \\
&= mz'_1-\tfrac{1}{2}z'_1(z'_1+1)+\tfrac{1}{2}n(n+1)-\tfrac{1}{2}k(k+1)-\tfrac{1}{2}m(m-1)  \\
&= \tfrac{1}{2}(n-k+m)(m-(n-k)-1)+\tfrac{1}{2}n(n+1)-\tfrac{1}{2}k(k+1)-\tfrac{1}{2}m(m-1) \\
&= k(n-k).
\end{align*}

\item Suppose that $m_1+m_2$ is odd.
Then ${\rm def}(\Lambda)=m_1-m_2$ is odd and negative, and
${\rm def}(\underline\theta(\Lambda))=m_2-m_1-1$ is even and non-negative.
By (\ref{0402}), we have
\begin{align*}
k &=\sum_{i=1}^{m_1}2a_i+\sum_{j=1}^{m_2}2b_j-\frac{1}{2}m(m-2)+\frac{1}{2}(m_2-m_1) \\
n &=2(\tau+m_1)+\sum_{i=1}^{m_1}2a_i+\sum_{j=1}^{m_2}2b_j-\frac{1}{2}(m+1)(m-1)+\frac{1}{2}(m_2-m_1-1)
\end{align*}
where $m=m_1+m_2$.
Then we have
\[
2(\tau+m_1)=n-k+m.
\]
Now by (\ref{0404}), we have
\begin{align*}
X^0_\Lambda&=\{2a_1,2a_2,\ldots,2a_{m_1}\},
& X^1_\Lambda&=\{2b_1+1,2b_2+1,\ldots,2b_{m_2}+1\}; \\
X^0_{\underline\theta(\Lambda)}&=\{2(\tau+m_1),2a_1,2a_2,\ldots,2a_{m_1}\},
& X^1_{\underline\theta(\Lambda)}&=\{2b_1+1,2b_2+1,\ldots,2b_{m_2}+1\}.
\end{align*}
The remaining proof is the same as in case (a).
\end{enumerate}

\item Suppose that $k+n$ is even.
The assumption that $\Lambda$ is $(n-k)$-admissible means that
\[
\underline\theta(\Lambda)=\binom{\tau+m_2,b_1,b_2,\ldots,b_{m_2}}{a_1,a_2,\ldots,a_{m_1}},\qquad
\Upsilon(\underline\theta(\Lambda))=\sqbinom{\tau,\nu_1,\nu_2,\ldots,\nu_{m_2}}{\mu_1,\mu_2,\ldots,\mu_{m_1}}.
\]
Now the proof is similar to that in case (1).
The only difference is that the element $\tau+m_1$ is replaced by $\tau+m_2$.
\end{enumerate}
\end{proof}

\begin{exam}
Consider the dual pair $(\rmU_k,\rmU_n)$ with $k=\frac{1}{2}d(d+1)$ and
$n=\frac{1}{2}(d+1)(d+2)$ where $d$ is a non-negative integer.
Let $\rho_d$ denote the unipotent cuspidal character of $\rmU_k$.
It is known that
\[
\deg_q(\rho_d)
=\frac{1}{24}(d-1)d(d+1)(3d+2).
\]
from \cite{lg} (9.5.1).
Note that $\rho_d=\rho_{\Lambda_d}$ where symbol $\Lambda_d$ is given in the proof of
Lemma~\ref{0410}.
It is easy to check that $\underline\theta(\rho_d)=\rho_{d+1}$ for the dual pair $(\rmU_k,\rmU_n)$,
and then
\begin{align*}
\deg_q(\rho_d)+k(n-k)
&= \frac{1}{24}(d-1)d(d+1)(3d+2)+\frac{1}{2}d(d+1)^2 \\
&= \frac{1}{24}d(d+1)(d+2)(3d+5) \\
&= \deg_q(\rho_{d+1}).
\end{align*}
\end{exam}

\section{Degree Differences in the Eta Correspondence}

\subsection{Lusztig correspondence for symplectic or orthogonal groups}\label{0510}
Let $\bfG$ be a symplectic group or an orthogonal group.
For $s\in (G^*)^0$, we define
\begin{align}\label{0515}
\begin{split}
\bfG^{(0)}=\bfG^{(0)}(s)
&=\prod_{\langle\lambda\rangle\subset\{\lambda_1,\ldots,\lambda_n\},\ \lambda\neq\pm 1}\bfG_{[\lambda]}(s);\\
\bfG^{(-)}=\bfG^{(-)}(s) &=\bfG_{[-1]}(s); \\
\bfG^{(+)}=\bfG^{(+)}(s) &=\begin{cases}
(\bfG_{[1]}(s))^*, & \text{if $\bfG$ is symplectic};\\
\bfG_{[1]}(s), & \text{otherwise},
\end{cases}
\end{split}
\end{align}
where $\bfG_{[\lambda]}(s)$ is given in \cite{amr} subsection 1.B
(see also \cite{pan-Lusztig-correspondence} subsection 2.2).
One can see that
\begin{itemize}
\item if $\bfG=\rmO^\epsilon_k$ where $k$ is even,
then $\bfG^{(-)}=\rmO^{\epsilon'}_{k^{(-)}}$, and $\bfG^{(+)}=\rmO^{\epsilon'\epsilon}_{k^{(+)}}$ where $k^{(-)},k^{(+)}$ are even;

\item if $\bfG=\Sp_k$,
then $\bfG^{(-)}=\rmO^{\epsilon'}_{k^{(-)}}$, and $\bfG^{(+)}=\Sp_{k^{(+)}}$ where $k^{(-)},k^{(+)}$ are even;

\item if $\bfG=\rmO_k$ where $k$ is odd,
then $\bfG^{(-)}=\Sp_{k^{(-)}}$, and $\bfG^{(+)}=\Sp_{k^{(+)}}$ where $k^{(-)},k^{(+)}$ are even,
\end{itemize}
for some $\epsilon'=+$ or $-$ and some $k^{(-)}, k^{(+)}$ all depending on $s$.
Note that $k^{(-)}+k^{(+)}\leq k$ for the first two cases above,
and $k^{(-)}+k^{(+)}\leq k-1$ if $\bfG=\rmO_k$ is an odd orthogonal group.
Then we have a (modified) \emph{Lusztig correspondence}
\begin{align*}
\Xi_s\colon \cale(G)_s &\rightarrow \cale(G^{(0)}(s)\times G^{(-)}(s)\times G^{(+)}(s))_1 \\
\rho &\mapsto \rho^{(0)}\otimes\rho_{\Lambda^{(-)}}\otimes\rho_{\Lambda^{(+)}}
\end{align*}
where $\Lambda^{(-)}\in\cals_{\bfG^{(-)}}$ and $\Lambda^{(+)}\in\cals_{\bfG^{(+)}}$.
The mapping $\Xi_s$ is two-to-one, more precisely $\Xi_s(\rho)=\Xi_s(\rho\cdot\sgn)$ when
$\bfG$ is an odd-orthogonal group;
and $\Xi_s$ is a bijection otherwise.

Let $(\bfG,\bfG')=(\rmO^\epsilon_k,\Sp_n)$ or $(\Sp_k,\rmO^\epsilon_n)$,
and let $\ell$ be a non-negative integer.
Now we define the $\ell$-admissibility of $\rho\in\cale(G)$ according the following cases:
\begin{enumerate}
\item[(I)] Suppose that $\bfG=\Sp_k$ where $k$ is even, $\ell$ is even.
Now $\bfG^{(-)}=\rmO^{\epsilon'}_{k^{(-)}}$  for some $\epsilon'=+$ or $-$,
and $\bfG^{(+)}=\Sp_{k^{(+)}}$.
Then $\rho\in\cale(G)$ is called \emph{$\ell$-admissible} for the dual pair
$(\bfG,\bfG')$ if
\begin{itemize}
\item $\rho_{\Lambda^{(-)}}$ is $\ell$-admissible for $(\bfG^{(-)},\Sp_{\rm even})$ and

\item $\rho_{\Lambda^{(+)}}$ is $\ell$-admissible for $(\bfG^{(+)},\rmO^{\epsilon'}_{\rm even})$
given in (\ref{0317}).
\end{itemize}

\item[(II)] Suppose that $\bfG=\Sp_k$ where $k$ is even, $\ell$ is odd.
Now $\bfG^{(-)}=\rmO^{\epsilon'}_{k^{(-)}}$  for some $\epsilon'=+$ or $-$,
and $\bfG^{(+)}=\Sp_{k^{(+)}}$.
Then $\rho\in\cale(G)$ is called \emph{$\ell$-admissible} for the dual pair
$(\bfG,\bfG')$ if
\begin{itemize}
\item $\rho_{\Lambda^{(-)}}$ is $(\ell-1)$-admissible for $(\bfG^{(-)},\Sp_{\rm even})$ and

\item $\rho_{\Lambda^{(+)}}$ is $(\ell-1)$-admissible for $(\bfG^{(+)},\rmO^{\epsilon'}_{\rm even})$.
\end{itemize}

\item[(III)] Suppose that $\bfG=\rmO^\epsilon_k$ where $k$ is even, $\ell$ is even.
Then both $\bfG^{(-)},\bfG^{(+)}$ are even-orthogonal groups.
Then $\rho\in\cale(G)$ is called \emph{$\ell$-admissible} for the dual pair
$(\bfG,\bfG')$ if
\begin{itemize}
\item $\rho_{\Lambda^{(-)}}$ is $\ell$-admissible for $(\bfG^{(-)},\Sp_{\rm even})$ and

\item $\rho_{\Lambda^{(+)}}$ is $\ell$-admissible for $(\bfG^{(+)},\Sp_{\rm even})$.
\end{itemize}

\item[(IV)] Suppose that $\bfG=\rmO_k$ where $k$ is odd, $\ell$ is odd.
Now $\bfG^{(-)}=\Sp_{k^{(-)}}$, and $\bfG^{(+)}=\Sp_{k^{(+)}}$.
Let $\rho\in\cale(G)$.
Because now $G\simeq\SO_k(q)\times\{\pm1\}$,
we can write $\rho\simeq\rho_1\otimes\rho_2$ where $\rho_1\in\cale(\SO_k(q))$
and $\rho_2\in\cale(\{\pm1\})$.
Let $\epsilon'=+$ of $\rho_2=\bf1$ and $\epsilon'=-$ if $\rho_2=\sgn$.
Then $\rho\in\cale(G)$ is called \emph{$\ell$-admissible} for the dual pair
$(\bfG,\bfG')$ if
\begin{itemize}
\item $\rho_{\Lambda^{(-)}}$ is $(\ell+1)$-admissible for $(\bfG^{(-)},\rmO^{\epsilon'}_{\rm even})$ and

\item $\rho_{\Lambda^{(+)}}$ is $(\ell-1)$-admissible for $(\bfG^{(+)},\rmO^{\epsilon'}_{\rm even})$.
\end{itemize}
\end{enumerate}
Note that the case when $\bfG=\rmO^\epsilon_k$ and $k+\ell$ is odd is undefined.
It is clear that if $\rho$ is $\ell$-admissible then it is also $\ell'$-admissible (if it is defined)
for any $\ell'\geq\ell$.

\begin{lem}\label{0516}
If the dual pair $(\bfG,\bfG')=(\Sp_k,\rmO^\epsilon_n)$ or $(\rmO^\epsilon_k,\Sp_n)$ is in stable range,
then every irreducible character $\rho\in\cale(G)$ is $(n-k)$-admissible.
\end{lem}
\begin{proof}
Suppose that $(\bfG,\bfG')$ is in stable range, $\rho\in\cale(G)_s$ for some $s$,
and write $\Xi_s(\rho)=\rho^{(0)}\otimes\rho_{\Lambda^{(-)}}\otimes\rho_{\Lambda^{(+)}}$.

\begin{enumerate}
\item Suppose that $(\bfG,\bfG')=(\Sp_k,\rmO^\epsilon_n)$ and $n$ is even.
Now $n-k$ is even and $k\leq n-k$.
Now both $k^{(\varepsilon)}\leq k\leq n-k$,
and so $\rho_{\Lambda^{(\varepsilon)}}$ is $(n-k)$-admissible by Lemma~\ref{0319} for
both $\varepsilon=+$ and $-$.

\item Suppose that $(\bfG,\bfG')=(\Sp_k,\rmO_n)$ and $n$ is odd.
Now $n-k$ is odd and $k\leq n-k-1$.
Now both $k^{(\varepsilon)}\leq k\leq n-k-1$,
and so $\rho_{\Lambda^{(\varepsilon)}}$ is $(n-k-1)$-admissible by Lemma~\ref{0319} for
both $\varepsilon=+$ and $-$.

\item Suppose that $(\bfG,\bfG')=(\rmO^\epsilon_k,\Sp_n)$ and $k$ is even.
Now $n-k$ is even and $k\leq n-k$.
Now both $k^{(\varepsilon)}\leq k\leq n-k$,
and so $\rho_{\Lambda^{(\varepsilon)}}$ is $(n-k)$-admissible by Lemma~\ref{0319} for
both $\varepsilon=+$ and $-$.

\item Suppose that $(\bfG,\bfG')=(\rmO_k,\Sp_n)$ and $k$ is odd.
Now $n-k$ is odd and $k\leq n-k$.
Now both $k^{(\varepsilon)}\leq k-1\leq n-k-1$,
and so $\rho_{\Lambda^{(\varepsilon)}}$ is $(n-k-1)$-admissible by Lemma~\ref{0319} for
both $\varepsilon=+$ and $-$.
\end{enumerate}
For all cases, we conclude that $\rho$ is $(n-k)$-admissible from the definition above.
\end{proof}

Now Lemma~\ref{0316} can be generalized to any irreducible characters:
\begin{lem}\label{0508}
Consider the dual pair $(\bfG,\bfG')=(\rmO^\epsilon_k,\Sp_n)$ or $(\Sp_k,\rmO^\epsilon_n)$ where $k\leq n$
and $\epsilon=+$ or $-$.
If $\rho\in\cale(G)$ is $(n-k)$-admissible,
then $\underline\theta(\rho)$ is defined and of\/ $\Theta$-rank $k$.
\end{lem}
\begin{proof}
Suppose that $\rho\in\cale(G)_s$ for some $s$ is $(n-k)$-admissible.
Write $\Xi_s(\rho)=\rho^{(0)}\otimes\rho_{\Lambda^{(-)}}\otimes\rho_{\Lambda^{(+)}}$.
Now we consider the following cases:
\begin{enumerate}
\item Suppose that $(\bfG,\bfG')=(\Sp_k,\rmO^\epsilon_n)$ where both $k,n$ are even and $k\leq n$.
By the definition above, both $\rho_{\Lambda^{(-)}},\rho_{\Lambda^{(+)}}$ are $(n-k)$-admissible.
Now $\bfG^{(-)}=\rmO^{\epsilon'}_{k^{(-)}}$ and $\bfG^{(+)}=\Sp_{k^{(+)}}$
for some $\epsilon'=+$ or $-$ depending on $s$.
Now $s$ is conjugate to an element of the form $(t,1_{k^{(+)}})$ where $t$ does not have
eigenvalue $1$.
Let $n^{(+)}$ be such that $n^{(+)}-k^{(+)}=n-k$ and so the unipotent character $\rho_{\Lambda^{(+)}}$
is $(n^{(+)}-k^{(+)})$-admissible.
Then $\underline\theta(\rho_{\Lambda^{(+)}})=\rho_{\underline\theta(\Lambda^{(+)})}$ is defined for the dual pair
$(\Sp_{k^{(+)}},\rmO^{\epsilon'\epsilon}_{n^{(+)}})$.
Hence $\underline\theta(\rho)$ is defined for the dual pair $(\bfG,\bfG')$ for some $s'\in G'^*$
via the following commutative diagram
\begin{equation}\label{0504}
\begin{CD}
\rho @> \underline\theta >> \underline\theta(\rho) \\
@V \Xi_s VV @VV \Xi_{s'} V \\
\rho^{(0)}\otimes\rho_{\Lambda^{(-)}}\otimes\rho_{\Lambda^{(+)}} @> \id\otimes\id\otimes\underline\theta >> \rho^{(0)}\otimes\rho_{\Lambda^{(-)}}\otimes\rho_{\underline\theta(\Lambda^{(+)})} \\
\end{CD}
\end{equation}
where $s'$ is an element of the form $(t,1_{n^{(+)}})$ so that $\bfG^{(0)}\simeq\bfG'^{(0)}$,
$\bfG^{(-)}\simeq\bfG'^{(-)}$ and $(\bfG^{(+)},\bfG'^{(+)})=(\Sp_{k^{(+)}},\rmO^{\epsilon'\epsilon}_{n^{(+)}})$.

As in the proof of Lemma~\ref{0316},
the $(n^{(+)}-k^{(+)})$-admissibility of $\rho_{\Lambda^{(+)}}$ implies that both unipotent characters
$\underline\theta(\rho_{\Lambda^{(+)}}),\underline\theta(\rho_{\Lambda^{(+)}})\sgn$ of
$\rmO^{\epsilon'\epsilon}_{n^{(+)}}(q)$ do not occur in the $\Theta$-correspondence for the dual pair
$(\Sp_{k'^{(+)}},\rmO^{\epsilon'\epsilon}_{n^{(+)}})$ for any $k'^{(+)}<k^{(+)}$.
Hence by \cite{pan-Lusztig-correspondence} theorem~6.9 and remark~6.10, we see that
both irreducible characters $\underline\theta(\rho),\underline\theta(\rho)\sgn$ of $\rmO^\epsilon_n(q)$
do not occur in the $\Theta$-correspondence for the dual pair $(\Sp_{k'},\rmO^\epsilon_n)$ for any $k'<k$.

We know that $\underline\theta(\rho)\chi_{\rmO^\epsilon_n}\in\cale(G')_{-s'}$,
and if we write 
\[
\Xi_{-s'}(\underline\theta(\rho)\chi_{\rmO^\epsilon_n})
=\rho'^{(0)}\otimes\rho_{\Lambda'^{(-)}}\otimes\rho_{\Lambda'^{(+)}},
\]
then we have $\Lambda'^{(-)}=\underline\theta(\Lambda^{(+)})$ and $\Lambda'^{(+)}=\Lambda^{(-)}$.
By the same argument as in the previous paragraph, the $(n^{(-)}-k^{(-)})$-admissibility of $\rho_{\Lambda^{(-)}}$
implies that both $\underline\theta(\rho)\chi_{\rmO^\epsilon_n},\underline\theta(\rho)\chi_{\rmO^\epsilon_n}\sgn$
of $\rmO^\epsilon_n(q)$ do not occur in the $\Theta$-correspondence for the dual pair $(\Sp_{k''},\rmO^\epsilon_n)$
for any $k''<k$.

\item Suppose that $(\bfG,\bfG')=(\Sp_k,\rmO_n)$ where $k$ is even, $n$ is odd, and $k\leq n$.
By definition, now both $\rho_{\Lambda^{(-)}},\rho_{\Lambda^{(+)}}$ are $(n-k-1)$-admissible.
Now $\bfG^{(-)}=\rmO^{\epsilon'}_{k^{(-)}}$ and $\bfG^{(+)}=\Sp_{k^{(+)}}$
for some $\epsilon'=+$ or $-$ depending on $s$.
Note that now both $k^{(-)},k^{(+)}$ are even.
Let $n^{(-)}$ be an even integer such that $n^{(-)}-k^{(-)}=n-k-1$ and so the unipotent character
$\rho_{\Lambda^{(-)}}$ is $(n^{(-)}-k^{(-)})$-admissible.
Then $\underline\theta(\rho_{\Lambda^{(-)}})=\rho_{\underline\theta(\Lambda^{(-)})}$ is defined for the dual pair
$(\rmO^{\epsilon'}_{k^{(-)}},\Sp_{n^{(-)}})$.
Hence $\underline\theta(\rho)$ is defined for the dual pair $(\Sp_k,\rmO_n)$
via the following commutative diagram
\begin{equation}\label{0505}
\begin{CD}
\rho @> \underline\theta >> \underline\theta(\rho) \\
@V \Xi_s VV @VV \iota\circ\Xi_{s'} V \\
\rho^{(0)}\otimes\rho_{\Lambda^{(-)}}\otimes\rho_{\Lambda^{(+)}} @> \id\otimes\underline\theta\otimes\id >> \rho^{(0)}\otimes\rho_{\underline\theta(\Lambda^{(-)})}\otimes\rho_{\Lambda^{(-)}} \\
\end{CD}
\end{equation}
where $s'\in G'^*$ is an element so that $\bfG^{(0)}\simeq\bfG'^{(0)}$,
$(\bfG^{(-)},\bfG'^{(+)})=(\rmO^{\epsilon'}_{k^{(-)}},\Sp_{n^{(+)}})$,
and $\bfG^{(+)}\simeq\bfG'^{(-)}$,
and $\iota(\rho^{(0)}\otimes\rho^{(-)}\otimes\rho^{(+)})=\rho^{(0)}\otimes\rho^{(+)}\otimes\rho^{(-)}$.
The remaining argument is similar to that in case (1).

\item Suppose that $(\bfG,\bfG')=(\rmO^\epsilon_k,\Sp_n)$ where both $k,n$ are even and $k\leq n$.
Now both $\rho_{\Lambda^{(-)}},\rho_{\Lambda^{(+)}}$ are $(n-k)$-admissible.
Moreover $\bfG^{(-)}=\rmO^{\epsilon'}_{k^{(-)}}$ and $\bfG^{(+)}=\rmO^{\epsilon'\epsilon}_{k^{(+)}}$
for some $\epsilon'=+$ or $-$ depending on $s$.
Let $n^{(+)}$ be such that $n^{(+)}-k^{(+)}=n-k$ and so the unipotent character $\rho_{\Lambda^{(+)}}$
is $(n^{(+)}-k^{(+)})$-admissible.
Then $\underline\theta(\rho_{\Lambda^{(+)}})=\rho_{\underline\theta(\Lambda^{(+)})}$ is defined for the dual pair
$(\Sp_{k^{(+)}},\rmO^{\epsilon'\epsilon}_{n^{(+)}})$.
Hence $\underline\theta(\rho)$ is defined for the dual pair $(\rmO^\epsilon_k,\Sp_n)$
via the same diagram (\ref{0504})
where $s'$ is an element so that $\bfG^{(0)}\simeq\bfG'^{(0)}$,
$\bfG^{(-)}\simeq\bfG'^{(-)}$ and $(\bfG^{(+)},\bfG'^{(+)})=(\rmO^{\epsilon'\epsilon}_{k^{(+)}},\Sp_{n^{(+)}})$.

As in the proof of Lemma~\ref{0316},
the $(n^{(+)}-k^{(+)})$-admissibility of $\rho_{\Lambda^{(+)}}$ implies that
$\underline\theta(\rho_{\Lambda^{(+)}})$ of
$\Sp_{n^{(+)}}(q)$ does not occur in the $\Theta$-correspondence for the dual pair
$(\rmO^{\epsilon''}_{k'^{(+)}},\Sp_{n^{(+)}})$ for any $\epsilon''=+$ or $-$ and any even $k'^{(+)}<k^{(+)}$.
Hence we see that the irreducible characters $\underline\theta(\rho)$ of $\Sp_n(q)$
does not occur in the $\Theta$-correspondence for the dual pair $(\rmO^{\epsilon''}_{k'},\Sp_n)$ 
for any $\epsilon''=+$ or $-$ and any even integer $k'$
such that $k'<k$.

Similarly, from diagram (\ref{0505}),
the $(n^{(-)}-k^{(-)})$-admissibility of $\rho_{\Lambda^{(-)}}$ implies that
the $\underline\theta(\rho)$ of $\Sp_n(q)$ does not occur in the $\Theta$-correspondence for the 
dual pair $(\rmO^{\epsilon''}_{k''},\Sp_n)$ for any $\epsilon''=+$ or $-$ and any odd integer $k''$
such that $k''<k$.

\item Suppose that $(\bfG,\bfG')=(\rmO^\epsilon_k,\Sp_n)$ where $k$ is odd, $n$ is even, and $k\leq n$.
By definition, now $\rho_{\Lambda^{(-)}}$ is $(n-k+1)$-admissible and $\rho_{\Lambda^{(+)}}$ is $(n-k-1)$-admissible.
Now $\bfG^{(-)}=\Sp_{k^{(-)}}$ and $\bfG^{(+)}=\Sp_{k^{(+)}}$.
Let $n^{(-)}$ be such that $n^{(-)}-k^{(-)}=n-k+1$ and so the unipotent character $\rho_{\Lambda^{(-)}}$
is $(n^{(-)}-k^{(-)})$-admissible.
Then $\underline\theta(\rho_{\Lambda^{(-)}})=\rho_{\underline\theta(\Lambda^{(-)})}$ is defined for the dual pair
$(\Sp_{k^{(-)}},\rmO^{\epsilon'}_{n^{(-)}})$.
Hence $\underline\theta(\rho)$ is defined for the dual pair $(\rmO^\epsilon_k,\Sp_n)$
via the commutative diagram
\[
\begin{CD}
\rho @> \underline\theta >> \underline\theta(\rho) \\
@V \iota\circ\Xi_s VV @VV \Xi_{s'} V \\
\rho^{(0)}\otimes\rho_{\Lambda^{(+)}}\otimes\rho_{\Lambda^{(-)}} @> \id\otimes\underline\theta\otimes\id >> \rho^{(0)}\otimes\rho_{\underline\theta(\Lambda^{(+)})}\otimes\rho_{\Lambda^{(+)}} \\
\end{CD}
\]
where $s'$ is an element so that $\bfG^{(0)}\simeq\bfG'^{(0)}$,
$(\bfG^{(+)},\bfG'^{(-)})=(\Sp_{k^{(+)}},\rmO^{\epsilon'}_{n^{(-)}})$,
and $\bfG^{(-)}\simeq\bfG'^{(+)}$.

As in the proof as above,
the $(n-k+1)$-admissibility of $\rho_{\Lambda^{(-)}}$ implies that
the irreducible characters $\underline\theta(\rho)$ of $\Sp_n(q)$
does not occur in the $\Theta$-correspondence for the dual pair $(\rmO^{\epsilon''}_{k'},\Sp_n)$ 
for any $\epsilon''=+$ or $-$ and any odd integer $k'$
such that $k'<k$.
Moreover, the $(n-k-1)$-admissibility of $\rho_{\Lambda^{(+)}}$ implies that the irreducible characters 
$\underline\theta(\rho)$ of $\Sp_n(q)$ does not occur in the $\Theta$-correspondence for the dual pair $(\rmO^{\epsilon''}_{k''},\Sp_n)$ for any $\epsilon''=+$ or $-$ and any even integer $k''$ such that $k''<k$.
\end{enumerate}
For all cases, we conclude that $\Theta\text{\rm -rk}(\underline\theta(\rho))=k$.
\end{proof}

\begin{prop}\label{0502}
Suppose that both $n,k$ are non-negative integers such that $k$ is even and $k\leq n$.
Then $\rho'\in\cale(\rmO^\epsilon_n(q))$ is of\/ $\Theta$-rank $k$ if and only if
there is a linear character $\chi$ of\/ $\rmO^\epsilon_n(q)$ such that
$\rho'\chi=\underline\theta(\rho)$ for some $(n-k)$-admissible character $\rho\in\cale(\Sp_k(q))$.
\end{prop}
\begin{proof}
Let $\rho'$ be an irreducible character of $\rmO^\epsilon_n(q)$.
Suppose that there is an $(n-k)$-admissible irreducible character $\rho$ of $\Sp_k(q)$
and a linear character $\chi$ of $\rmO^\epsilon_n(q)$ such that $\underline\theta(\rho)=\rho'\chi$.
By Lemma~\ref{0508}, we know that $\Theta\text{\rm -rk}(\rho'\chi)=k$
and hence $\Theta\text{\rm -rk}(\rho')=k$.

Conversely, suppose that $\rho'\in\cale(\rmO^\epsilon_n(q))$ such that $\Theta\text{\rm -rk}(\rho')=k$.
From \cite{pan-theta-rank} we know that there exists a linear character $\chi$ of $\rmO^\epsilon_n(q)$ such 
that $\rho'\chi=\underline\theta(\rho)$ for some $\rho\in\cale(\Sp_k(q))$.
Now we need to show that $\rho$ is $(n-k)$-admissible.
Suppose that $\rho\in\cale(\Sp_k(q))_s$ for some $s$ and
write $\Xi_s(\rho)=\rho^{(0)}\otimes\rho_{\Lambda^{(-)}}\otimes\rho_{\Lambda^{(+)}}$
where $\Lambda^{(-)}\in\cals_{\rmO^{\epsilon'}_{k^{(-)}}}$ and $\Lambda^{(+)}\in\cals_{\Sp_{k^{(+)}}}$
for some $\epsilon'=+$ or $-$ depending on $s$.
\begin{enumerate}
\item Suppose that $n$ is even.
If $\rho_{\Lambda^{(+)}}$ is not $(n-k)$-admissible,
then by the proof of Proposition~\ref{0302} we see that $\rho'\chi$ or $\rho'\chi\sgn$ have 
$\Theta$-rank less than $k$.
Similarly, if $\rho_{\Lambda^{(-)}}$ is not $(n-k)$-admissible,
then $\rho'\chi\chi_{\rmO^\epsilon_n}$ or $\rho'\chi\chi_{\rmO^\epsilon_n}\sgn$ will have
$\Theta$-rank less than $k$.

\item Suppose that $n$ is odd.
If $\rho_{\Lambda^{(-)}}$ is not $(n-k-1)$-admissible,
then $\rho'\chi$ or $\rho'\chi\sgn$ will have $\Theta$-rank less than $k$.
Similarly, if $\rho_{\Lambda^{(+)}}$ is not $(n-k-1)$-admissible,
then $\rho'\chi\chi_{\rmO^\epsilon_n}$ or $\rho'\chi\chi_{\rmO^\epsilon_n}\sgn$ will have
$\Theta$-rank less than $k$.
\end{enumerate}
Therefore, we conclude that $\rho$ must be $(n-k)$-admissible.
\end{proof}

\begin{prop}\label{0509}
Suppose that both $n,k$ are non-negative integers such that $n$ is even and $k\leq n$.
Then $\rho'\in\cale(\Sp_n(q))$ is of\/ $\Theta$-rank $k$ is and only if
$\rho'=\underline\theta(\rho)$ for some $(n-k)$-admissible character $\rho\in\cale(\rmO^\epsilon_k(q))$
and some $\epsilon=+$ or $-$.
\end{prop}
\begin{proof}
Let $\rho'$ be an irreducible character of $\Sp_n(q)$.
Suppose that there is an $(n-k)$-admissible irreducible character $\rho$ of $\rmO^\epsilon_k(q)$
for $\epsilon=+$ or $-$ such that $\underline\theta(\rho)=\rho'$.
By Lemma~\ref{0508}, we know that $\Theta\text{\rm -rk}(\rho')=k$.

Conversely, suppose that $\rho'\in\cale(\Sp_n(q))$ has $\Theta$-rank $k$.
Then there is an orthogonal group $\rmO^\epsilon_k(q)$ and $\rho\in\cale(\rmO^\epsilon_k(q))$
such that $\underline\theta(\rho)=\rho'$.
Suppose that $\rho\in\cale(\rmO^\epsilon_k(q))_s$ for some $s$ and
write $\Xi_s(\rho)=\rho^{(0)}\otimes\rho_{\Lambda^{(-)}}\otimes\rho_{\Lambda^{(+)}}$.
\begin{enumerate}
\item Suppose that $k$ is even.
Then $\Lambda^{(-)}\in\cals_{\rmO^{\epsilon'}_{k^{(-)}}}$ and
$\Lambda^{(+)}\in\cals_{\rmO^{\epsilon'\epsilon}_{k^{(+)}}}$ for some $\epsilon'=+$ or $-$ depending on $s$.
If $\rho_{\Lambda^{(+)}}$ is not $(n-k)$-admissible,
then $\rho'$ will occur in the $\Theta$-correspondence for the dual pair
$(\rmO^\epsilon_{k'},\Sp_n)$ for some even integer $k'<k$.
Similarly, if $\rho_{\Lambda^{(-)}}$ is not $(n-k)$-admissible,
then $\rho'$ will occur in the $\Theta$-correspondence for the dual pair
$(\rmO^\epsilon_{k'},\Sp_n)$ for some odd integer $k'<k$.

\item Suppose that $k$ is odd.
Then $\Lambda^{(-)}\in\cals_{\Sp_{k^{(-)}}}$ and $\Lambda^{(+)}\in\cals_{\Sp_{k^{(+)}}}$.
If $\rho_{\Lambda^{(+)}}$ is not $(n-k-1)$-admissible,
then $\rho'$ will occur in the $\Theta$-correspondence for the dual pair
$(\rmO^\epsilon_{k'},\Sp_n)$ for some odd integer $k'<k$.
Similarly, if $\rho_{\Lambda^{(-)}}$ is not $(n-k+1)$-admissible,
then $\rho'$ will occur in the $\Theta$-correspondence for the dual pair
$(\rmO^\epsilon_{k'},\Sp_n)$ for some even integer $k'<k$.
\end{enumerate}
Therefore, $\rho$ is $(n-k)$-admissible and so the proposition is proved.
\end{proof}

Now we generalize Proposition~\ref{0303} to general irreducible characters.

\begin{prop}\label{0512}
Consider the dual pair $(\bfG,\bfG')=(\rmO^\epsilon_k,\Sp_n)$ or $(\Sp_k,\rmO^\epsilon_n)$ where $k\leq n$.
If $\rho\in\cale(G)$ is $(n-k)$-admissible,
then
\[
\deg_q(\underline\theta(\rho))=\deg_q(\rho)+
\begin{cases}
\frac{1}{2}k(n-k+1), & \text{if\/ $(\bfG,\bfG')=(\rmO^\epsilon_k,\Sp_n)$};\\
\frac{1}{2}k(n-k-1), & \text{if\/ $(\bfG,\bfG')=(\Sp_k,\rmO^\epsilon_n)$}.
\end{cases}
\]
\end{prop}
\begin{proof}
For $\rho\in\cale(G)_s$ for some $s$,
we write $\Xi_s(\rho)=\rho^{(0)}\otimes\rho_{\Lambda^{(-)}}\otimes\rho_{\Lambda^{(+)}}$.
From \cite{lusztig-book} 4.23 (see also \cite{DM} remark 13.24),
it is known that
\begin{equation}\label{0503}
\rho(1)=
\frac{|G|_{p'}}{|G^{(0)}|_{p'}|G^{(-)}|_{p'}|G^{(+)}|_{p'}}
\rho^{(0)}(1)\rho_{\Lambda^{(-)}}(1)\rho_{\Lambda^{(+)}}(1).
\end{equation}

\begin{enumerate}
\item
Suppose that $(\bfG,\bfG')=(\rmO^\epsilon_k,\Sp_n)$ or $(\Sp_k,\rmO^\epsilon_n)$
where both $k,n$ are even.
Now we have a commutative diagram
\[
\begin{CD}
\rho @> \underline\theta >> \rho' \\
@V \Xi_s VV @VV \Xi_{s'} V \\
\rho^{(0)}\otimes\rho_{\Lambda^{(-)}}\otimes\rho_{\Lambda^{(+)}} @> \id\otimes\id\otimes\underline\theta >> \rho'^{(0)}\otimes\rho_{\Lambda'^{(-)}}\otimes\rho_{\Lambda'^{(+)}} \\
\end{CD}
\]
We know that $\bfG^{(0)}\simeq\bfG'^{(0)}$, $\rho^{(0)}=\rho'^{(0)}$
and $\bfG^{(-)}\simeq\bfG'^{(-)}$, $\rho_{\Lambda^{(-)}}=\rho_{\Lambda'^{(-)}}$.
Therefore by (\ref{0503}),
we have
\begin{multline*}
\deg_q(\underline\theta(\rho))-\deg_q(\rho)
=\deg_q(\underline\theta(\rho_{\Lambda^{(+)}}))-\deg_q(\rho_{\Lambda^{(+)}}) \\
+\deg_q(|G'|_{p'})-\deg_q(|G|_{p'})
-\deg_q(|G'^{(+)}|_{p'})+\deg_q(|G^{(+)}|_{p'}).
\end{multline*}

\begin{enumerate}
\item
Suppose that $(\bfG,\bfG')=(\rmO^\epsilon_k,\Sp_n)$.
Then $\bfG^{(+)}=\rmO^\epsilon_{k-2l}$ and $\bfG'^{(+)}=\Sp_{n-2l}$ for some non-negative integer $l$.
Hence by Proposition~\ref{0303}, we have
\begin{align*}
\deg_q(\underline\theta(\rho_{\Lambda^{(+)}}))-\deg_q(\rho_{\Lambda^{(+)}})
&= \tfrac{1}{2}(k-2l)(n-2l-(k-2l)+1), \\
\deg_q(|G'|_{p'})-\deg_q(|G|_{p'})
&= \tfrac{1}{4}n(n+2)-\tfrac{1}{4}k^2, \\
-\deg_q(|G'^{(+)}|_{p'})+\deg_q(|G^{(+)}|_{p'})
&=-\tfrac{1}{4}(n-2l)(n-2l+21)+\tfrac{1}{4}(k-2l)^2.
\end{align*}
Therefore,
\begin{align*}
& \deg_q(\underline\theta(\rho))-\deg_q(\rho) \\
&= \tfrac{1}{2}(k-2l)(n-k+1)+\tfrac{1}{4}n(n+2)-\tfrac{1}{4}k^2
-\tfrac{1}{4}(n-2l)(n-2l+2)+\tfrac{1}{4}(k-2l)^2 \\
&=\tfrac{1}{2}k(n-k+1).
\end{align*}

\item
Suppose that $(\bfG,\bfG')=(\Sp_k,\rmO^\epsilon_n)$.
Then $\bfG^{(+)}=\Sp_{k-2l}$ and $\bfG'^{(+)}=\rmO^\epsilon_{n-2l}$ for some non-negative integer $l$.
Hence by Proposition~\ref{0303}, we have
\begin{align*}
\deg_q(\underline\theta(\rho_{\Lambda^{(+)}}))-\deg_q(\rho_{\Lambda^{(+)}})
&= \tfrac{1}{2}(k-2l)(n-2l-(k-2l)-1), \\
\deg_q(|G'|_{p'})-\deg_q(|G|_{p'})
&= \tfrac{1}{4}n^2-\tfrac{1}{4}k(k+2), \\
-\deg_q(|G'^{(+)}|_{p'})+\deg_q(|G^{(+)}|_{p'})
&=-\tfrac{1}{4}(n-2l)^2+\tfrac{1}{4}(k-2l)(k-2l+2).
\end{align*}
Therefore
\begin{align*}
& \deg_q(\underline\theta(\rho))-\deg_q(\rho) \\
&= \tfrac{1}{2}(k-2l)(n-k-1)+\tfrac{1}{4}n^2-\tfrac{1}{4}k(k+2)
-\tfrac{1}{4}(n-2l)^2+\tfrac{1}{4}(k-2l)(k-2l+21) \\
&= \tfrac{1}{2}k(n-k-1).
\end{align*}
\end{enumerate}

\item Suppose that $(\bfG,\bfG')=(\rmO_k,\Sp_n)$ where $k$ is odd and $n$ even;
or $(\Sp_k,\rmO_n)$ where $k$ is even and $n$ is odd.
\begin{enumerate}
\item Suppose that $(\bfG,\bfG')=(\rmO_k,\Sp_n)$ where $k$ is odd and $n$ is even.
Now we have a commutative diagram
\[
\begin{CD}
\rho @> \underline\theta >> \rho' \\
@V \iota\circ\Xi_s VV @VV \Xi_{s'} V \\
\rho^{(0)}\otimes\rho_{\Lambda^{(+)}}\otimes\rho_{\Lambda^{(-)}} @> \id\otimes\underline\theta\otimes\id >> \rho'^{(0)}\otimes\rho_{\Lambda'^{(-)}}\otimes\rho_{\Lambda'^{(+)}} \\
\end{CD}
\]
We know that $\bfG^{(0)}\simeq\bfG'^{(0)}$, $\rho^{(0)}=\rho'^{(0)}$
and $\bfG^{(-)}\simeq\bfG'^{(+)}$, $\rho_{\Lambda^{(-)}}=\rho_{\Lambda'^{(+)}}$.
Therefore by (\ref{0503}),
we have
\begin{multline*}
\deg_q(\underline\theta(\rho))-\deg_q(\rho)
=\deg_q(\underline\theta(\rho_{\Lambda^{(+)}}))-\deg_q(\rho_{\Lambda^{(+)}}) \\
+\deg_q(|G'|_{p'})-\deg_q(|G|_{p'})
-\deg_q(|G'^{(-)}|_{p'})+\deg_q(|G^{(+)}|_{p'}).
\end{multline*}
Then $\bfG^{(+)}=\Sp_{k-2l-1}$ and $\bfG'^{(-)}=\rmO^\epsilon_{n-2l}$ for some non-negative integer $l$.
Hence by Proposition~\ref{0303}, we have
\begin{align*}
\deg_q(\underline\theta(\rho_{\Lambda^{(-)}}))-\deg_q(\rho_{\Lambda^{(-)}})
&= \tfrac{1}{2}(k-2l-1)(n-2l-(k-2l-1)-1), \\
\deg_q(|G'|_{p'})-\deg_q(|G|_{p'})
&= \tfrac{1}{4}n(n+2)-\tfrac{1}{4}(k-1)(k+1), \\
-\deg_q(|G'^{(-)}|_{p'})+\deg_q(|G^{(+)}|_{p'})
&=-\tfrac{1}{4}(n-2l)^2+\tfrac{1}{4}(k-1-2l)(k+1-2l).
\end{align*}
Therefore
\begin{align*}
\deg_q(\underline\theta(\rho))-\deg_q(\rho)
&= \tfrac{1}{2}(k-2l-1)(n-k)+\tfrac{1}{4}n(n+2)-\tfrac{1}{4}(k-1)(k+1) \\
&\qquad\qquad\qquad\qquad -\tfrac{1}{4}(n-2l)^2+\tfrac{1}{4}(k-1-2l)(k+1-2l) \\
&= \tfrac{1}{2}k(n-k+1).
\end{align*}

\item Suppose that $(\bfG,\bfG')=(\Sp_k,\rmO_n)$ where $k$ is even and $n$ is odd.
Now we have a commutative diagram
\[
\begin{CD}
\rho @> \underline\theta >> \rho' \\
@V \Xi_s VV @VV \iota\circ\Xi_{s'} V \\
\rho^{(0)}\otimes\rho_{\Lambda^{(-)}}\otimes\rho_{\Lambda^{(+)}} @> \id\otimes\underline\theta\otimes\id >> \rho'^{(0)}\otimes\rho_{\Lambda'^{(+)}}\otimes\rho_{\Lambda'^{(-)}} \\
\end{CD}
\]
We know that $\bfG^{(0)}\simeq\bfG'^{(0)}$, $\rho^{(0)}=\rho'^{(0)}$
and $\bfG^{(+)}\simeq\bfG'^{(-)}$, $\rho_{\Lambda^{(+)}}=\rho_{\Lambda'^{(-)}}$.
Therefore by (\ref{0503}),
we have
\begin{multline*}
\deg_q(\underline\theta(\rho))-\deg_q(\rho)
=\deg_q(\underline\theta(\rho_{\Lambda^{(-)}}))-\deg_q(\rho_{\Lambda^{(-)}}) \\
+\deg_q(|G'|_{p'})-\deg_q(|G|_{p'})
-\deg_q(|G'^{(+)}|_{p'})+\deg_q(|G^{(-)}|_{p'}).
\end{multline*}

Then $\bfG^{(-)}=\rmO^\epsilon_{k-2l}$ and $\bfG'^{(+)}=\Sp_{n-2l-1}$ for some non-negative integer $l$.
Hence by Proposition~\ref{0303}, we have
\begin{align*}
\deg_q(\underline\theta(\rho_{\Lambda^{(-)}}))-\deg_q(\rho_{\Lambda^{(-)}})
&= \tfrac{1}{2}(k-2l)(n-2l-1-(k-2l)+1), \\
\deg_q(|G'|_{p'})-\deg_q(|G|_{p'})
&= \tfrac{1}{4}(n-1)(n+1)-\tfrac{1}{4}k(k+2), \\
-\deg_q(|G'^{(+)}|_{p'})+\deg_q(|G^{(-)}|_{p'})
&=-\tfrac{1}{4}(n-1-2l)(n+1-2l)+\tfrac{1}{4}(k-2l)^2.
\end{align*}
Therefore
\begin{align*}
\deg_q(\underline\theta(\rho))-\deg_q(\rho)
&= \tfrac{1}{2}(k-2l)(n-k)+\tfrac{1}{4}(n-1)(n+1)-\tfrac{1}{4}k(k+2) \\
&\qquad\qquad\qquad\quad -\tfrac{1}{4}(n-1-2l)(n+1-2l)+\tfrac{1}{4}(k-2l)^2 \\
&=\tfrac{1}{2}k(n-k-1).
\end{align*}
\end{enumerate}
\end{enumerate}
\end{proof}

\subsection{Lusztig correspondence for unitary groups}
Let $\bfG=\rmU_k$.
For $s\in G^*=\rmU_k(q)$, we can write
\begin{equation}\label{0506}
C_{G^*}(s)=\prod_{i=1}^r\prod_{j=1}^{t_i}\GL_{k_{ij}}^{(-1)^i}(q^i)
\end{equation}
for some non-negative integers $r,t_i,k_{ij}$ such that $\sum_{i=1}^r\sum_{j=1}^{t_i}ik_{ij}=k$ where
$\GL^{+1}_{k_{ij}}=\GL^{+}_{k_{ij}}:=\GL_{k_{ij}}$
and $\GL^{-1}_{k_{ij}}=\GL^{-}_{k_{ij}}:=\rmU_{k_{ij}}$.
For $i=1$, we can let $t_1=q+1$ and each $j=1,\ldots,t_1$ corresponds an
eigenvalue $\lambda_j\in\overline\bff_q$ of $s$ such that $\lambda_j^{q+1}=1$.
Let $\rho\in\cale(G)_s$ and for the Lusztig correspondence
$\grL_s\colon\cale(G)_1\rightarrow\cale(C_{G^*}(s))_1$ we can write
\begin{equation}\label{0308}
\grL_s(\rho)=\bigotimes_{i=1}^r\bigotimes_{j=1}^{t_i}\rho^{(ij)}
\end{equation}
for some $\rho^{(ij)}\in\cale(\GL_{k_{ij}}^{(-1)^i}(q^i))_1$.
Note that $\rho^{(1j)}$ is a unipotent character of $\rmU_{k_{1j}}(q)$.

Now for a non-negative integer $\ell$,
an irreducible character $\rho\in\cale(\rmU_k(q))$ is called \emph{$\ell$-admissible}
if each $\rho^{(1j)}$ is $\ell$-admissible for $j=1,\ldots,t_1$.
It is clear that if $\rho$ is $\ell$-admissible,
then $\rho$ is also $\ell'$-admissible for any $\ell'\geq\ell$,
and $\rho\chi$ is also $\ell$-admissible for any linear character $\chi$ of $\rmU_k(q)$.

\begin{lem}\label{0517}
If the dual pair $(\bfG,\bfG')=(\rmU_k,\rmU_n)$ is in stable range,
then every irreducible character $\rho\in\cale(G)$ is $(n-k)$-admissible.
\end{lem}
\begin{proof}
Suppose that the dual pair $(\bfG,\bfG')=(\rmU_k,\rmU_n)$ is in stable range,
and let $\rho\in\cale(G)_s$ for some $s$.
Write $\grL_s(\rho)=\bigotimes_{i=1}^r\bigotimes_{j=1}^{t_i}\rho^{(ij)}$ as above.
Let $k_{ij}$ be given as in (\ref{0506}).
Then $k_{1j}\leq k\leq n-k$ and hence $\rho^{(1j)}$ is $(n-k)$-admissible by Lemma~\ref{0407}
for each $j=1,\ldots,t_1$.
Hence $\rho$ is $(n-k)$-admissible by the definition above.
\end{proof}

\begin{lem}\label{0513}
Consider the dual pair $(\bfG,\bfG')=(\rmU_k,\rmU_n)$ where $k\leq n$.
If $\rho\in\cale(G)$ is $(n-k)$-admissible,
then $\underline\theta(\rho)$ is defined and of\/ $\Theta$-rank $k$.
\end{lem}
\begin{proof}
Let $\rho$ be an $(n-k)$-admissible character  of $\rmU_k(q)$.
Suppose that $\rho\in\cale(G)_s$ for some $s$.
Let $\bfG^{(1)}=U_{k_{1,1}}$ be the component in (\ref{0506}) corresponding the eigenvalue $1$,
and let $\bfG^{(0)}$ be the product of the other components.
So we have $C_{G^*}(s)=G^{(0)}\times G^{(1)}$ and
$\grL_s(\rho)=\rho^{(0)}\otimes\rho^{(1)}$ where $\rho^{(1)}=\rho^{(1,1)}$
and $\rho^{(0)}$ is the tensor product of the components $\rho^{(ij)}$ with $(i,j)\neq(1,1)$.

Let $n_{1,1}$ be such that $n_{1,1}-k_{1,1}=n-k$,
so the unipotent character $\rho^{(1,1)}$ is $(n_{1,1}-k_{1,1})$-admissible.
Then $\underline\theta(\rho^{(1)})$ is defined for the dual pair
$(\rmU_{k_{1,1}},\rmU_{n_{1,1}})$ by Lemma~\ref{0414}.
Hence $\underline\theta(\rho)$ is defined for the dual pair $(\rmU_k,\rmU_n)$
via the following commutative diagram

\begin{equation}\label{0507}
\begin{CD}
\rho @> \underline\theta >> \rho' \\
@V \grL_s VV @VV \grL_{s'} V \\
\rho^{(0)}\otimes\rho^{(1)} @> \id\otimes\underline\theta >> \rho^{(0)}\otimes\underline\theta(\rho^{(1)}) \\
\end{CD}
\end{equation}
where $s'=(s,1_{n-k})\in\rmU_k(q)\times\rmU_{n-k}(q)\subset\rmU_n(q)$.

Let $\chi$ be a linear character of $\rmU_n(q)$ and suppose that $\underline\theta(\rho)\chi$
is in $\cale(\rmU_n(q))_{s''}$ for some $s''$.
Note that if we write
\[
\grL_{s'}(\underline\theta(\rho))=\bigotimes_{i=1}^{r'}\bigotimes_{j=1}^{t'_i}\rho'^{(ij)},\qquad
\grL_{s''}(\underline\theta(\rho)\chi)=\bigotimes_{i=1}^{r''}\bigotimes_{j=1}^{t''_i}\rho''^{(ij)},
\]
then we know that two ordered set of the components
\[
\{\rho'^{(1,1)},\rho'^{(1,2)},\ldots,\rho'^{(1t'_1)}\}\quad\text{and}\quad
\{\rho''^{(1,1)},\rho''^{(1,2)},\ldots,\rho''^{(1t''_1)}\}
\]
are different by a permutation.
Now the assumption that each $\rho^{(1j)}$ is $(n-k)$-admissible means that
$\underline\theta(\rho)\chi$ of $\rmU_n(q)$ does not occurs
in the $\Theta$-correspondence for the dual pair $(\rmU_{k'},\rmU_n)$ for any $k'<k$.
This implies that $\Theta\text{\rm -rk}(\underline\theta(\rho))=k$.
\end{proof}

\begin{prop}\label{0514}
Let $n,k$ be non-negative integers such that $k\leq n$.
Then $\rho'\in\cale(\rmU_n(q))$ is of\/ $\Theta$-rank $k$ if and only if
there is a linear character $\chi$ of\/ $\rmU_n(q)$ such that
$\rho'\chi=\underline\theta(\rho)$ for some $(n-k)$-admissible character $\rho\in\cale(\rmU_k(q))$.
\end{prop}
\begin{proof}
The proof is similar to those of Proposition~\ref{0502} or Proposition~\ref{0509}.
\end{proof}

Now we generalize Proposition~\ref{0401} to general irreducible characters.

\begin{prop}\label{0511}
Consider the dual pair $(\bfG,\bfG')=(\rmU_k,\rmU_n)$ where $k\leq n$.
If an irreducible character $\rho\in\cale(\rmU_k(q))$ is $(n-k)$-admissible,
then
\[
\deg_q(\underline\theta(\rho))=\deg_q(\rho)+k(n-k).
\]
\end{prop}
\begin{proof}
Let $\rho$ be an $(n-k)$-admissible irreducible character  of $\rmU_k(q)$.
Suppose that $\rho\in\cale(G)_s$ for some $s$ and write $\grL_s(\rho)=\rho^{(0)}\otimes\rho^{(1)}$.
Let $\rho'=\underline\theta(\rho)\in\cale(G')_{s'}$ for some $s'$ and write
$\grL_{s'}(\rho')=\rho'^{(0)}\otimes\rho'^{(1)}$ as above.
Now we have a commutative diagram
\[
\begin{CD}
\rho @> \underline\theta >> \rho' \\
@V \grL_s VV @VV \grL_{s'} V \\
\rho^{(0)}\otimes\rho^{(1)} @> \id\otimes\underline\theta >> \rho'^{(0)}\otimes\rho'^{(1)} \\
\end{CD}
\]
We know that
\begin{align*}
\deg_q(\underline\theta(\rho))
&= |G'|_{p'}+\deg_q(\rho'^{(0)})+\deg_q(\rho'^{(1)})-|G'^{(0)}|_{p'}-|G'^{(1)}|_{p'}, \\
\deg_q(\rho)
&= |G|_{p'}+\deg_q(\rho^{(0)})+\deg_q(\rho^{(1)})-|G^{(0)}|_{p'}-|G^{(1)}|_{p'}.
\end{align*}
Now we have $\bfG^{(1)}=\rmU_{k-l}$ and $\bfG'^{(1)}=\rmU_{n-l}$ for some non-negative integer $l\leq k$,
$G^{(0)}\simeq G'^{(0)}$, $\rho^{(0)}=\rho'^{(0)}$.
Then by Proposition~\ref{0401}, we have
\begin{align*}
\deg_q(\rho'^{(1)})-\deg_q(\rho^{(1)})
&= (k-l)(n-l-(k-l))=(k-l)(n-k), \\
\deg_q(|G'|_{p'})-\deg_q(|G|_{p'})
&= \tfrac{1}{2}n(n+1)-\tfrac{1}{2}k(k+1), \\
-\deg_q(|G'^{(1)}|_{p'})+\deg_q(|G^{(1)}|_{p'})
&=-\tfrac{1}{2}(n-l)(n-l+1)+\tfrac{1}{2}(k-l)(k-l+1).
\end{align*}
Therefore
\begin{align*}
&\deg_q(\underline\theta(\rho))-\deg_q(\rho) \\
&=(k-l)(n-k)+\tfrac{1}{2}n(n+1)-\tfrac{1}{2}k(k+1)
-\tfrac{1}{2}(n-l)(n-l+1)+\tfrac{1}{2}(k-l)(k-l+1) \\
&=k(n-k).
\end{align*}
\end{proof}

\subsection{Proofs of the main results}

\begin{proof}[Proof of Theorem~\ref{0101}]
Part (i) is Lemma~\ref{0508} and Lemma~\ref{0513}.
Part (ii) is Proposition~\ref{0502}, Proposition~\ref{0509} and Proposition~\ref{0514}.
\end{proof}

\begin{proof}[Proof of Theorem~\ref{0501}]
The theorem is Proposition~\ref{0512} and Proposition~\ref{0511}.
\end{proof}

\bibliography{refer}
\bibliographystyle{amsalpha}

\end{document}